\documentclass[11pt,leqno]{amsart}
\usepackage[margin=1.1in]{geometry}
\usepackage{amssymb,amsmath,amsfonts,amsthm}
\usepackage{tikz}
\usepackage{needspace}
\usepackage[hidelinks]{hyperref}
\newtheorem{theorem}{Theorem}[section]
\newtheorem{lemma}[theorem]{Lemma}
\newtheorem{proposition}[theorem]{Proposition}
\newtheorem{corollary}[theorem]{Corollary}
\theoremstyle{definition}
\newtheorem{example}[theorem]{Example}
\theoremstyle{remark}
\newtheorem{remark}[theorem]{Remark}
\numberwithin{equation}{section}
\DeclareMathOperator{\Hol}{Hol}
\DeclareMathOperator{\Tr}{Tr}
\DeclareMathOperator{\rank}{rank}
\newcommand{\R}{\mathbb R}
\newcommand{\C}{\mathbb C}
\newcommand{\F}{\mathcal F}
\title{Inverse spectra and Wilson traces of $SU(2)$ connections on finite graphs}
\author{Haoxuan Cheng}
\address{School of Mathematical Sciences, Fudan University,
Shanghai 200433, China}
\email{hxcheng25@m.fudan.edu.cn}
\date{}
\begin{document}
\raggedbottom
\begin{abstract}
We study a fixed $SU(2)$ connection on a known finite labelled simple graph
through its full Laplacian spectra with independently variable positive edge
weights. This family is equivalent to all vertex-simple cycle traces.
For connected graphs, it determines every vertex gauge class exactly for
trees and graphs with a single cyclic block that is a cycle, a theta graph,
or a subdivision of a simple three-vertex-connected cubic graph.
For arbitrary graphs, we describe the full fibers, including degenerate
connections, and give criteria for generic uniqueness and Wilson
observability under restricted edge controls.
\end{abstract}
\maketitle

\section{Introduction}
Let $G=(V,E)$ be a known finite connected labelled simple graph;
thus there are no loops or parallel edges, and all vertex and edge labels
are fixed. The group $SU(2)$ consists of the determinant-one unitary
$2\times2$ matrices. A connection $U$ assigns a matrix $U_{yx}\in SU(2)$
to transport a vector in $\C^2$ from $x$ to $y$, with
$U_{xy}=U_{yx}^{-1}$ for the reverse edge. Choosing one orientation per
edge identifies the connection space with $SU(2)^E$. Along a path
$x_0,\ldots,x_k$, the transport is
$U_{x_kx_{k-1}}\cdots U_{x_1x_0}$: the first traversed edge acts first,
on the right. For a closed path $C$, this product is its holonomy,
denoted by $\Hol_U(C)$.

A vertex gauge is a choice $g_x\in SU(2)$ at every vertex, acting by
$U_{yx}\mapsto g_yU_{yx}g_x^{-1}$. Two connections are gauge equivalent
if one such choice transforms all their edge matrices at once. We seek
to recover this equivalence class, called the vertex gauge class.
This is also the unit quaternion gain-graph model, in which vertex gauge
is called switching \cite{BelardoEtAl,CavaleriDonno}.
Quaternionic spectral representations and constructions of cospectral gain
graphs are treated in \cite{CavaleriDonnoSpessato}.

For positive real edge weights $w=(w_e)_{e\in E}$, write $w_{xy}$ for
the weight of $\{x,y\}$. With vertex measure one, the connection
Laplacian acts on fields $f:V\to\C^2$ by
\[
 (L_w(U)f)(x)=\sum_{y:\{x,y\}\in E}w_{xy}\bigl(f(x)-U_{xy}f(y)\bigr).
\]
A vertex gauge conjugates this operator by the corresponding block
unitary matrix and therefore preserves its spectrum. We keep $U$ fixed,
vary all positive weights independently, and record the full spectrum,
with multiplicities, at each known weight vector. Does this
weight-labelled family determine the vertex gauge class of $U$?

Connection Laplacians also occur in vector diffusion maps \cite{SingerWu}
and group synchronization \cite{BandeiraSingerSpielman}, where pairwise
edge transformations are input data. The inverse problem here asks what
can be recovered from the weight-labelled spectra alone.

The spectral problem has an exact reduction to loop data.
For a closed path $C$, its Wilson trace is
$\tau_C(U)=\tfrac12\operatorname{tr}\Hol_U(C)$, where $\operatorname{tr}$
is the ordinary $2\times2$ matrix trace. This real number is invariant
under reversal of $C$, change of starting vertex, and vertex gauge.
These loop observables originate in lattice gauge theory \cite{Wilson}.
Proposition~\ref{prop:equivalence}
shows that the entire weight-labelled spectral family contains exactly
the Wilson traces of the vertex-simple cycles, which visit each of their
vertices only once. The forward reduction isolates a cycle by varying
weights; the converse uses Kenyon's self-dual determinant expansion
\cite[Section 9.1, Theorem 8]{Kenyon}. Related determinant formulas for Laplacians of quaternion unit gain graphs
are given in \cite{KyrcheiTreisterPelykh}.
Wilson traces over all closed walks
are known to separate gauge classes \cite[Theorem 3.1]{Levy}.
The restriction to simple cycles is thus the issue to resolve here.

A pair of cycles meeting at one vertex illustrates the loss of information
(Figure~\ref{fig:graphclasses}(a)).
Each simple cycle lies in one of the two cycles, so the observations determine
the conjugacy class of each loop holonomy separately. One may conjugate one
holonomy while keeping the other fixed and preserve all these traces.
Yet a gauge at the common vertex conjugates both simultaneously, and the
two resulting pairs need not be equivalent. The way cycles meet can
therefore obstruct recovery even when every positive weight has been tested.

\subsection{Graphs with recovery at every connection}
The main theorem gives an exact graph criterion for unique recovery.
Separation of gauge classes means that any two connections with the same
observations are gauge equivalent. A cut vertex is a vertex whose deletion
increases the number of connected components. A graph is biconnected if
it is connected and has no cut vertex; a cyclic block is a maximal
biconnected subgraph containing a cycle.

A cubic graph has degree three at every vertex.
A theta graph consists of three internally vertex-disjoint paths with
the same two endpoints. A subdivision replaces edges by paths with new
internal vertices. Three-vertex-connectivity requires at least four vertices
and connectivity after deleting any set of at most two vertices.

\Needspace{10\baselineskip}
\begin{samepage}
\begin{theorem}\label{thm:allconnections}
For a finite connected labelled simple graph $G$, the following are equivalent:
\begin{enumerate}
\item All simple-cycle traces separate every pair of vertex gauge classes.
\item All independently variable positive edge-weight spectra separate every
pair of vertex gauge classes.
\item $G$ is a tree, or it has exactly one cyclic block, which is a cycle,
a theta graph, or a subdivision of a finite simple three-vertex-connected
cubic graph.
\end{enumerate}
\end{theorem}
\end{samepage}

\begin{figure}[htbp]
\centering
\begin{tikzpicture}[x=.95cm,y=.95cm,font=\small,line width=.65pt,
  vertex/.style={circle,fill=black,inner sep=1.6pt,outer sep=0pt},
  subdiv/.style={circle,draw,fill=white,inner sep=2pt,outer sep=0pt}]
  \foreach \n/\x/\y in {a0/0/0,a1/-.95/.75,a2/-.95/-.75,
                          a3/.95/.75,a4/.95/-.75}
    \node[vertex] (\n) at (\x,\y) {};
  \draw (a0)--(a1)--(a2)--(a0)--(a3)--(a4)--(a0);
  \node[below=3pt] at (a0) {$v$};
  \node at (0,-1.35) {(a)};
  \foreach \n/\x/\y in {b0/3.1/.9,b1/2.25/.28,b2/2.57/-.73,
                          b3/3.63/-.73,b4/3.95/.28}
    \node[vertex] (\n) at (\x,\y) {};
  \draw (b0)--(b1)--(b2)--(b3)--(b4)--(b0);
  \node at (3.1,-1.35) {(b)};
  \foreach \n/\x/\y in {c0/5.25/0,c1/7.15/0,c2/6.2/.85,
                          c3/6.2/0,c4/6.2/-.85}
    \node[vertex] (\n) at (\x,\y) {};
  \draw (c0)--(c2)--(c1) (c0)--(c3)--(c1) (c0)--(c4)--(c1);
  \node at (6.2,-1.35) {(c)};
  \foreach \n/\x/\y in {d0/9.1/1,d1/8.1/-.9,d2/10.1/-.9,
                          d3/9.1/-.3,t0/10.9/-.25,t1/11.6/.25,t2/11.7/-.75}
    \node[vertex] (\n) at (\x,\y) {};
  \node[subdiv] (ds) at (9.1,-.9) {};
  \draw (d0)--(d1)--(ds)--(d2)--(d0)
        (d0)--(d3) (d1)--(d3)--(d2);
  \draw[dashed] (d2)--(t0)--(t1) (t0)--(t2);
  \node at (9.9,-1.35) {(d)};
\end{tikzpicture}
\caption{Graphs in the main theorem. (a) Two cycles share a cut vertex $v$ and admit connections without unique recovery; (b) a cycle; (c) a theta graph; (d) a subdivided $K_4$ with an attached tree. The open vertex subdivides an edge; dashed edges form the attached tree. Panels (b)--(d) illustrate Theorem~\ref{thm:allconnections}; $K_4$ is just one representative of the three-connected cubic cores. All drawings are schematic: lengths and angles have no metric meaning.}
\label{fig:graphclasses}
\end{figure}
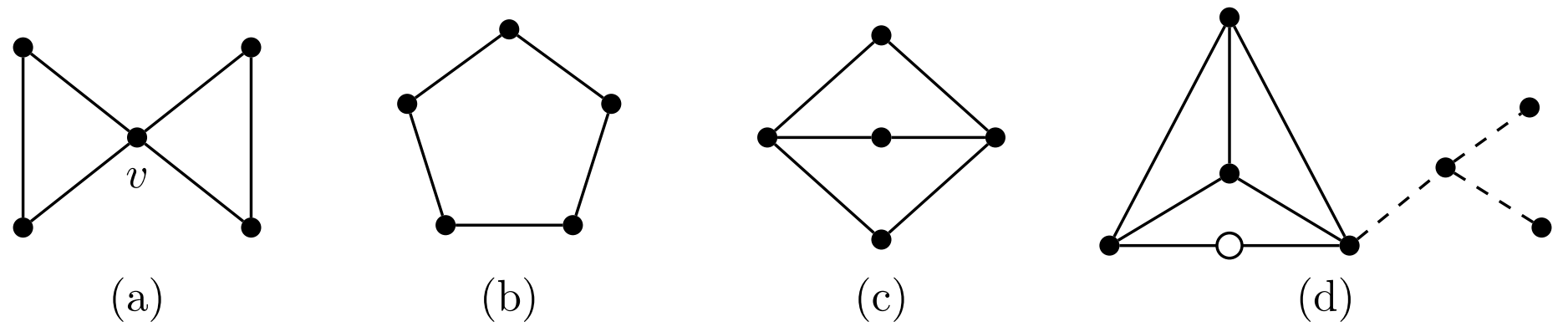

A tree has one gauge class. In condition (3), a graph with cycles is equivalently obtained by attaching trees to the listed cyclic block, each tree meeting that block at one vertex. Figure~\ref{fig:graphclasses}(b)--(d) illustrates the three cyclic structures; panel (d) also shows subdivision and tree attachment, starting from $K_4$,
the complete graph on four vertices. Section~\ref{sec:allconnections} proves the theorem.

The proof must transfer rigidity through changes of the cubic graph at
every connection, including the exceptional ones. One difficulty is that a
simple cycle after a local change can become a closed walk through two vertex-disjoint old
cycles joined by an edge. We show that these cycles and the specified edge
lie in a subdivision of a triangular prism, the graph consisting of two
disjoint triangles joined at corresponding vertices by three edges. Its simple-cycle data determine its
restricted gauge class, hence also the trace of that joined walk. This
allows the data comparison to pass through the local change. When the trace equations lose rank, we determine which remaining choices
differ only by rotations fixing the old transports.
The calculations use unit quaternion trace coordinates
\cite{Watson} and the graph gauge-quotient viewpoint \cite{Ramadas}.

\subsection{The ambiguity when recovery is not unique}
We also determine the full fiber of the observations: the set of all gauge
classes producing the same data. For probabilistic statements, Haar
measure on $SU(2)$ is its normalized translation-invariant probability
measure. A property is \emph{generic} if it holds on a gauge-invariant open
dense subset of $SU(2)^E$ of full product Haar measure.
The full-measure assertions below use the zero-set property of nontrivial
real-analytic functions in local coordinates \cite{Mityagin}.

On every cyclic block the fiber is finite and can be constructed
recursively. At a generic connection it has one or two classes, with a
graph-minor criterion for the singleton case. A minor is obtained by
deleting vertices or edges and contracting edges
(Section~\ref{sec:generic}). For a general connected graph, all
continuous ambiguity comes from placing compatible block connections
relative to one another, modulo a common rotation
(Theorem~\ref{thm:everyfiber}). The two-cycle example is one instance of
this relative freedom. We describe how compatible classes can merge as the
data vary and when continuous data can be realized by a continuously varying gauge
class. Such a realization can fail even when every individual datum is
realizable (Section~\ref{sec:continuousdata}).

\subsection{Finite measurements and restricted controls}
The use of an entire weight family raises a further question: how much of
it must actually be measured? Section~\ref{sec:finite} proves that finitely
many full spectra preserve every fiber. It also selects a number of cycle
traces linear in the cycle rank $\beta=|E|-|V|+1$, the number of
independent cycles, and gives certified small-noise recovery of
the entire fiber at generic connections. The target of reconstruction is
therefore a class when uniqueness holds and the complete compatible set
when it does not.

If only some edge weights can be varied, the reduced family may lose Wilson
traces as well. Section~\ref{sec:arbitrarycontrols} gives an exact decision
procedure for whether it still recovers all the original simple-cycle
traces, for any control set and fixed positive weights algebraic over $\mathbb Q$.
Combining this decision with the fiber description distinguishes information
lost in measurement from ambiguity already present with full control.
Operator extensions and higher-rank obstructions are treated in
Section~\ref{sec:extensions}; the longer estimates and coefficient
calculations are in the appendices.

\section{Spectral data and quaternion coordinates}\label{sec:model}
Keep the connection, path order, gauge action and Laplacian convention
introduced above. Write $n=|V|$, $m=|E|$ and $\beta=m-n+1$; for a
connected subgraph $H$, its cycle rank is $|E(H)|-|V(H)|+1$.
The compact gauge quotient of this connection space is
$X_G=SU(2)^E/SU(2)^V$, and $[U]\in X_G$ denotes the class of $U$.
The graph labels stay fixed under this quotient. We write $Q^*$ for the
conjugate transpose and $I_d$ for the $d$-dimensional identity matrix.

The scalar degree diagonal of $L_w(U)$ changes with $w$. Write
$\lambda_1(w,U)\le\cdots\le\lambda_{2n}(w,U)$ for its eigenvalues,
listed with multiplicity. The all-weight spectral observation is the map
\[
 \Sigma_G(U):(0,\infty)^E\longrightarrow\R^{2n},\qquad
 w\longmapsto\bigl(\lambda_1(w,U),\ldots,\lambda_{2n}(w,U)\bigr).
\]
Thus $\Sigma_G(U)=\Sigma_G(V)$ means equality of spectra at the same
weight vector for every $w>0$. The parameter $w$ is part of the data.

Let $\mathcal C(G)$ be the finite set of vertex-simple cycles, with neither
orientation nor starting vertex specified. For $Q\in SU(2)$ put
$\tau(Q)=\tfrac12\operatorname{tr}Q$, where $\operatorname{tr}$ is the
ordinary $2\times2$ matrix trace, and set
\[
 \tau_C(U)=\tau(\Hol_U(C)),\qquad
 W_G:X_G\longrightarrow[-1,1]^{\mathcal C(G)},\quad
 [U]\longmapsto(\tau_C(U))_{C\in\mathcal C(G)}.
\]
The $SU(2)$ trace is real and unchanged by inversion or conjugation, so
these definitions are independent of the cycle choices and of the gauge.
We write $W_G(U)$ for $W_G([U])$. For a fixed connection $U$, its full
observation fiber is
\[
 \F_U:=W_G^{-1}(W_G(U))
       =\{[U']\in X_G:W_G(U')=W_G(U)\}.
\]
A fiber therefore consists of all gauge classes compatible with the data.
The proposition below identifies it with the fiber of $\Sigma_G$.
Throughout, $\Tr$ denotes the trace on the full $2n$-dimensional vertex
space, and $\|v\|_\infty=\max_j|v_j|$ for a finite real coordinate vector.

\begin{proposition}\label{prop:equivalence}
For all connections $U,V$, $\Sigma_G(U)=\Sigma_G(V)$ if and only if
$W_G(U)=W_G(V)$.
\end{proposition}
\begin{proof}
For a cycle $C$ of length $\ell$, set its edge weights to one and all
other weights to $t>0$; denote this weight vector by $w_C(t)$.
Let $I$ also denote the trivial connection, with every transport equal to
$I_2$. Define the scalar spectral-moment polynomial
$F_U(t)=\Tr L_{w_C(t)}(U)^\ell$. It is a polynomial of degree at most $\ell$, so positive measurements determine
$F_U(0)$. At zero the graph is $C$ and isolated vertices. A walk winding
around $C$ in $\ell$ steps uses every cycle edge and no diagonal step;
there are $\ell$ starting points in each orientation. A closed walk with
zero winding reduces to backtracking pairs, whose transports cancel;
inserting scalar diagonal steps does not change this. Thus only the
$2\ell$ winding walks depend on the connection. Each contributes
$(-1)^\ell\operatorname{tr}\Hol_U(C)=(-1)^\ell2\tau_C(U)$ to the full
trace, giving the factor $4\ell$ below. Therefore
\begin{equation}\label{eq:isolation}
 F_U(0)-F_I(0)=(-1)^\ell4\ell\bigl(\tau_C(U)-1\bigr).
\end{equation}
For the converse, first take $z\in\R$. For a $2\times2$ complex block
$B$, let $\widetilde B=\operatorname{adj}B$ be its adjugate. A block
matrix $M=(M_{xy})$ is self-dual when $M_{yx}=\widetilde{M_{xy}}$;
its diagonal blocks must then be scalar. Since
$\widetilde U=U^{-1}$ for $U\in SU(2)$ and
$\widetilde{wU}=wU^{-1}$ for real $w$, the matrix $zI-L_w(U)$ has this
structure. This is quaternionic self-duality in the $2\times2$ complex
representation.
Kenyon's permutation expansion
\cite[Section 9.1, Theorem 8]{Kenyon} gives its quaternion determinant
using simple-cycle traces for cycles of length at least three, scalar
diagonal entries for one-cycles and scalar edge norms for two-cycles.
The complex determinant is its square. Equal simple-cycle data thus give
equal characteristic polynomials. These polynomial identities extend to all $z\in\C$, including at
singular matrices and zero edge blocks.
\end{proof}

\subsection{Quaternion coordinates}
We identify $SU(2)$ with the unit sphere $S^3\subset\R^4$ of quaternions.
Here $S^d=\{v\in\R^{d+1}:|v|=1\}$, with $|\cdot|$ the Euclidean norm.
Using the complex imaginary unit $i$, the Pauli matrices are
\[
 \sigma_x=\begin{pmatrix}0&1\\1&0\end{pmatrix},\quad
 \sigma_y=\begin{pmatrix}0&-i\\i&0\end{pmatrix},\quad
 \sigma_z=\begin{pmatrix}1&0\\0&-1\end{pmatrix}.
\]
The quaternion $(a,u)\in\R\times\R^3$ represents the matrix
$aI_2+i(u_1\sigma_x+u_2\sigma_y+u_3\sigma_z)$; its scalar and imaginary
parts are $\operatorname{Re}(a,u)=a$ and $\operatorname{Im}(a,u)=u$.
In this basis the product is
\[
 (a,u)(b,v)=(ab-u\cdot v,\ av+bu-u\times v),\qquad
 \operatorname{Re}(PQ^{-1})=\langle P,Q\rangle\quad(P,Q\in S^3).
\]
The dot product, cross product and $\langle\cdot,\cdot\rangle$ are the
usual Euclidean operations in $\R^3$ and $\R^4$, respectively.
Trace invariants under simultaneous matrix conjugation are treated by
Procesi \cite{Procesi}; Goldman \cite{Goldman} develops trace coordinates
for determinant-one complex $2\times2$ matrices.
A Gram matrix records pairwise inner products; $\det(u,v,z)$ means the
determinant with the three indicated vectors as columns.
We write $O(d)$ for the real orthogonal group and $SO(d)$ for its
subgroup of determinant one. A rotation acts on quaternions by
$R\cdot(a,u)=(a,Ru)$.
Conjugation induces $SO(3)$ on imaginary parts. Independent left and right
unit multiplications induce all of $SO(4)$: they preserve orientation,
and their six independent infinitesimal generators span its Lie algebra.
Their kernel is the common pair of signs. A generating-tree gauge sets
tree edges to identity. The edges outside the tree, called chords, then
carry $\beta$ arbitrary unit quaternions modulo one common conjugation. Conditional on the tree transports, these chord
variables are independent Haar variables, meaning samples from the normalized
translation-invariant probability measure on $SU(2)$. The effective
stabilizer of a chord tuple consists of the rotations in $SO(3)$ fixing
every quaternion in that tuple. A rotation fixing a nonzero imaginary
vector is a rotation about its axis. Hence the effective stabilizer is
$SO(3)$ if all imaginary parts vanish, $SO(2)$ if their span is a line,
and trivial if their span has dimension at least two.

We call a connection central if its based holonomy subgroup lies in
$\{\pm1\}$, and commuting if that subgroup is abelian. In a tree gauge
these conditions concern the chord tuple; they do not require the original
edge matrices to commute.

\section{Proof of all-connection rigidity}\label{sec:allconnections}

The proof has a positive and a negative part. For the positive part,
we prove rigidity of two small cores and show that it survives the local
operations generating simple three-connected cubic graphs. For the negative
part, we first exclude graphs with multiple cyclic blocks, then identify
the cuts producing reflection ambiguities inside a single cyclic block. Throughout this section we write
$\mathsf P(H)$ for separation of every gauge class by the simple cycles of $H$.
For the cycle-word calculations below, put $T_{xy}=U_{yx}^{-1}$ and multiply
these variables in path order. This gives the inverse of the original path
transport and preserves every closed-path trace. The gauge action is
$T_{xy}\mapsto g_xT_{xy}g_y^{-1}$.
In the explicit quaternion examples, choose perpendicular imaginary units
$i,j$ and put $k=ij$. This fixes the Hamilton convention for those examples.

\subsection{Two rigid cores}
\begin{lemma}\label{lem:k4}
The seven simple-cycle traces of $K_4$ determine every connection up to gauge.
\end{lemma}
\begin{proof}
Gauge the three edges incident to vertex $0$ to identity and write
$A=U_{21}$, $B=U_{31}$, $C=U_{32}$. The observations are half-traces of
\[
 A,\ B,\ C,\ BA^{-1},\ C^{-1}B,\ CA,\ B^{-1}CA.
\]
For $A=(a,u)$, $B=(b,v)$, $C=(c,z)$, the first six give all scalar
parts and the imaginary Gram matrix:
\[
 u\cdot v=\tau(BA^{-1})-ab,\quad
 v\cdot z=\tau(C^{-1}B)-bc,\quad u\cdot z=ac-\tau(CA).
\]
Here $\tau(Q)=\operatorname{tr}Q/2$. The last gives
\[
 \det(u,v,z)=abc-b(u\cdot z)+c(u\cdot v)+a(v\cdot z)-\tau(B^{-1}CA).
\]
The diagonal Gram entries are also known, since $|u|^2=1-a^2$,
$|v|^2=1-b^2$ and $|z|^2=1-c^2$. Two tuples with the same Gram matrix
are related by an isometry of their spans, which extends to an element
of $O(3)$. At rank three the displayed determinant forces this element
to have determinant one. At lower rank, if necessary, compose it with
a reflection fixing the target span pointwise. The resulting map lies
in $SO(3)$ and still matches all three vectors. Quaternion conjugation
realizes this rotation, so the connections are gauge equivalent in every case.
\end{proof}

\begin{lemma}\label{lem:prismall}
The triangular prism and every subdivision of it have property $\mathsf P$.
\end{lemma}
\begin{proof}
Use triangles $012,345$ and links $03,14,25$, as labelled in
Figure~\ref{fig:prism}(a). Gauge the tree
$01,02,03,14,25$ to one and set
$A=T_{12},B=T_{34},C=T_{35},D=T_{45}$.
The fourteen simple cycles give the half-traces of
\[
\begin{gathered}
A,B,C,D,\quad AC^{-1},B^{-1}C,DC^{-1},A^{-1}B^{-1},D^{-1}B^{-1},AD^{-1},\\
AD^{-1}B^{-1},A^{-1}DC^{-1},AC^{-1}B,BDC^{-1}.
\end{gathered}
\]
The single variables give the scalar parts; unit norms then give the
diagonal Gram entries. If $X=(x_0,x)$ and $Y=(y_0,y)$, then
$\tau(X^\epsilon Y^\delta)=x_0y_0-\epsilon\delta\,x\cdot y$ for
$\epsilon,\delta\in\{1,-1\}$. The six two-variable words therefore give
the six distinct imaginary inner products. For a third variable
$Z=(z_0,z)$, our multiplication convention gives
\[
 \tau(XYZ)=x_0y_0z_0-z_0(x\cdot y)-x_0(y\cdot z)
                 -y_0(x\cdot z)+\det(x,y,z).
\]
Replacing a variable by its inverse changes the sign of its imaginary
part. Thus the last four words determine, respectively, the volumes for
the triples $(A,B,D)$, $(A,C,D)$, $(A,B,C)$ and $(B,C,D)$, with
coefficients $-1,-1,1,1$. All other terms in these four equations are
already known. Hence all triple volumes, as well as the full Gram matrix,
are determined.
These are complete simultaneous $SO(3)$ invariants: a Gram isometry is
proper in rank three by the volumes, and can be made proper on the
orthogonal complement in lower rank. Subdivision preserves path transports
and adds only internal vertex gauges.
\end{proof}

\subsection{Transferring data through local moves}
A local change of a cubic core need not preserve which closed walks are
simple. We therefore need to recover one additional kind of trace from
the original simple-cycle data: the trace obtained by traversing one cycle,
crossing a joining edge, traversing a disjoint cycle, and returning along
that edge. The prism provides this information while retaining the
specified joining edge.

\begin{figure}[htbp]
\centering
\begin{tikzpicture}[x=1cm,y=1cm,font=\small,line width=.65pt,
  vertex/.style={circle,fill=black,inner sep=1.7pt,outer sep=0pt},
  subdiv/.style={circle,draw,fill=white,inner sep=2pt,outer sep=0pt}]
  \node[vertex,label=above:$0$] (p0) at (0,1.65) {};
  \node[vertex,label=below left:$1$] (p1) at (-1.65,-1) {};
  \node[vertex,label=below right:$2$] (p2) at (1.65,-1) {};
  \node[vertex,label=above:$3$] (p3) at (0,.6) {};
  \node[vertex,label=below left:$4$] (p4) at (-.6,-.35) {};
  \node[vertex,label=below right:$5$] (p5) at (.6,-.35) {};
  \draw (p0)--(p1)--(p2)--(p0) (p3)--(p4)--(p5)--(p3)
        (p0)--(p3) (p1)--(p4) (p2)--(p5);
  \node at (0,-1.55) {(a)};
  \begin{scope}[xshift=6.3cm]
    \node[vertex,label=below:$u$] (u) at (-1.05,0) {};
    \node[vertex] (l1) at (-2.4,.9) {};
    \node[vertex] (l2) at (-2.4,-.9) {};
    \node[vertex,label=below:$v$] (v) at (1.05,0) {};
    \node[vertex] (r1) at (2.4,.9) {};
    \node[vertex] (r2) at (2.4,-.9) {};
    \draw (u)--(l1)--(l2)--(u) (v)--(r1)--(r2)--(v);
    \draw[line width=1.25pt] (u)--node[above] {$e$}(v);
    \node at (-1.95,0) {$C_1$};
    \node at (1.95,0) {$C_2$};
    \coordinate (top) at (0,1.55);
    \coordinate (bottom) at (0,-1.55);
    \draw[dashed] (l1) to[out=45,in=180] (top)
      to[out=0,in=135] (r1);
    \draw[dashed] (l2) to[out=-45,in=180] (bottom)
      to[out=0,in=-135] (r2);
    \node[subdiv] at (top) {};
    \node[subdiv] at (bottom) {};
    \node at (0,-2) {(b)};
  \end{scope}
\end{tikzpicture}
\caption{The prism and completion of a specified dumbbell. (a) Triangles $012,345$ and links $03,14,25$, with the labels used in Lemma~\ref{lem:prismall}. (b) Lemma~\ref{lem:dumbbellall} retains both vertex-disjoint cycles $C_1,C_2$ and the specified edge $e=uv$ (thick), and supplies two further paths (dashed), forming a prism subdivision. Open vertices indicate possible internal vertices; path lengths may vary.}
\label{fig:prism}
\end{figure}
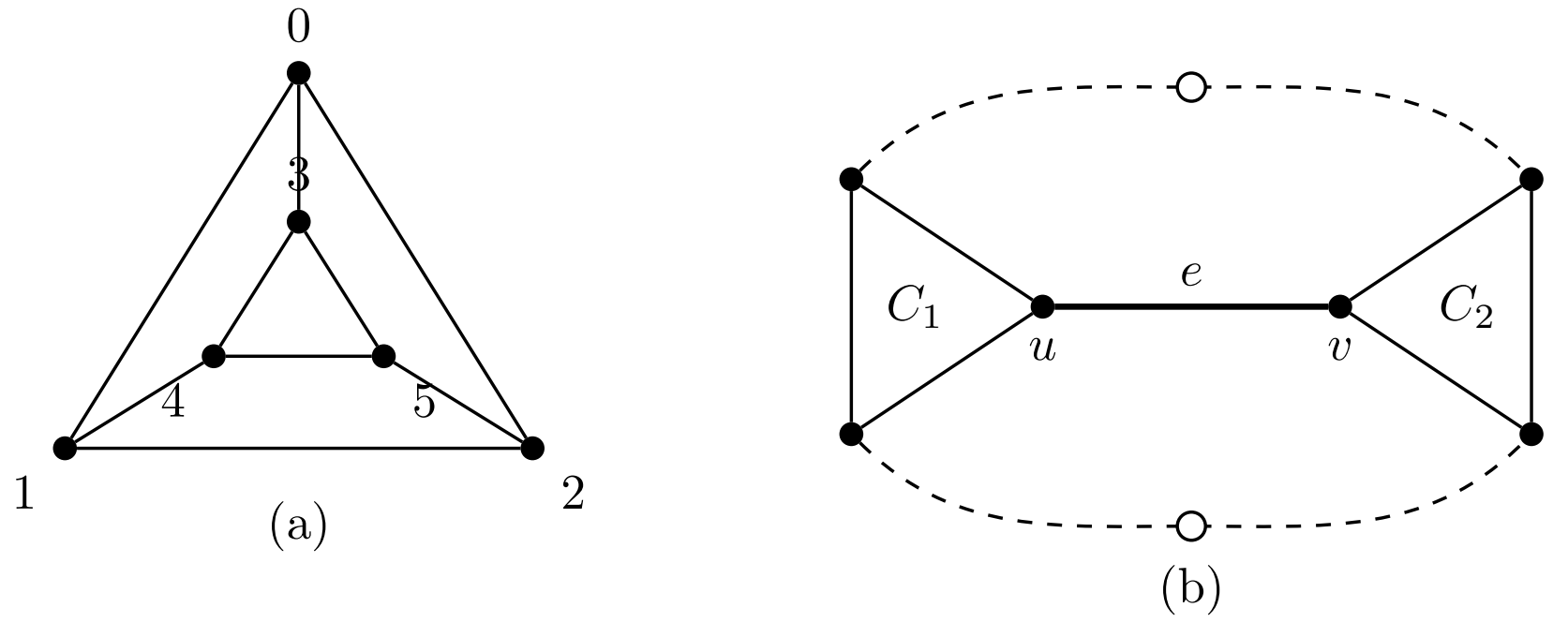

\begin{lemma}\label{lem:dumbbellall}
In a simple three-connected cubic graph, two vertex-disjoint cycles joined
by a specified edge lie in a prism subdivision containing that entire
edge and both cycles. Their joined dumbbell trace is determined by simple cycles.
\end{lemma}
\begin{proof}
Let $e=uv$ join $C_1,C_2$. Three-vertex-connectivity implies
three-edge-connectivity. For $w\ne u,v$, suppose $G-w-e$ had $k\ge2$
components. Each is a proper vertex set with at least three boundary
edges in $G$. Summing their boundaries counts edges incident to $w$
at most once and $e$ at most twice, giving at most $\deg w+2=5$,
contrary to $3k\ge6$. Thus $G-w-e$ is connected.

Now put $A=V(C_1)\setminus\{u\}$ and $B=V(C_2)\setminus\{v\}$.
If $G-\{u,v\}$ had no two vertex-disjoint $A$--$B$ paths, the set version
of Menger's theorem would give an $A$--$B$ separator of size at most one.
Write it as $\{w\}$; if it is empty, delete any additional vertex other
than $u,v$. Since the cycles have length at least three, both $A-w$ and
$B-w$ remain nonempty. In $G-\{u,v,w\}$ no component meets both sets.
Cubicity makes the two neighbors of $u$ other than $v$ lie on $C_1$,
and the two neighbors of $v$ other than $u$ lie on $C_2$. Restoring $u,v$
without $e$ can therefore join components meeting $A-w$ only among
themselves, and likewise on the $B-w$ side; it cannot join the two sides.
This contradicts connectivity of $G-w-e$.

Take the two vertex-disjoint connecting paths. On each, stop at its first
contact with $C_2$ and start at its last preceding contact with $C_1$.
The resulting paths have interiors outside both cycles. Their endpoints
on either cycle are distinct, and none is $u$ or $v$, because the paths
were chosen in $G-\{u,v\}$. Together with $e$ they are three disjoint
links with three distinct endpoints on each cycle. The arcs of each
cycle between these endpoints form the two subdivided triangles of a
prism. Thus the union is the required prism subdivision, shown
schematically in Figure~\ref{fig:prism}(b). Its simple cycles are simple
cycles of $G$, so Lemma~\ref{lem:prismall} recovers its restricted gauge
class and hence the joined dumbbell trace.
\end{proof}

\begin{figure}[htbp]
\centering
\begin{tikzpicture}[x=.9cm,y=.9cm,font=\small,line width=.65pt,>=stealth,
  vertex/.style={circle,fill=black,inner sep=1.8pt,outer sep=0pt},
  port/.style={circle,draw,fill=white,inner sep=1.6pt,outer sep=0pt}]
  \foreach \idx/\s in {0/0,1/4.2,2/8.4} {
    \begin{scope}[shift={(\s,0)}]
      \node[port,label=above:$1$] (a\idx) at (-1.2,.85) {};
      \node[port,label=below:$2$] (b\idx) at (-1.2,-.85) {};
      \node[port,label=above:$3$] (c\idx) at (1.2,.85) {};
      \node[port,label=below:$4$] (d\idx) at (1.2,-.85) {};
    \end{scope}
  }
  \node[vertex] (oldl) at (-.45,0) {};
  \node[vertex] (oldr) at (.45,0) {};
  \draw (a0)--(oldl)--(b0) (c0)--(oldr)--(d0);
  \draw[line width=1.15pt] (oldl)--node[below=2pt] {$e$}(oldr);
  \node[vertex] (contract) at (4.2,0) {};
  \draw (a1)--(contract)--(b1) (c1)--(contract)--(d1);
  \node[vertex] (newt) at (8.4,.45) {};
  \node[vertex] (newb) at (8.4,-.45) {};
  \draw (a2)--(newt)--(c2) (b2)--(newb)--(d2);
  \draw[line width=1.15pt] (newt)--node[right=2pt] {$e'$}(newb);
  \draw[->] (1.65,0)--(2.55,0);
  \draw[->] (5.85,0)--(6.75,0);
  \node at (4.2,-1.5) {(a)};
  \begin{scope}[shift={(1,-3.5)}]
    \node[vertex,label=above:$v$] (v) at (0,0) {};
    \node[port,label=left:$r_1$] (r1) at (-1.4,1) {};
    \node[port,label=left:$r_2$] (r2) at (-1.4,-1) {};
    \node[port,label=right:$r_3$] (r3) at (1.4,0) {};
    \draw (r1)--(v)--(r2) (v)--(r3);
  \end{scope}
  \begin{scope}[shift={(7.4,-3.5)}]
    \node[vertex,label=above:$x_1$] (x1) at (-.6,.65) {};
    \node[vertex,label=below:$x_2$] (x2) at (-.6,-.65) {};
    \node[vertex,label=below right:$x_3$] (x3) at (.8,0) {};
    \node[port,label=left:$r_1$] (s1) at (-1.6,1) {};
    \node[port,label=left:$r_2$] (s2) at (-1.6,-1) {};
    \node[port,label=right:$r_3$] (s3) at (1.8,0) {};
    \draw (s1)--(x1)--(x2)--(s2) (x1)--(x3)--(s3);
    \draw[dashed] (x2)--(x3);
  \end{scope}
  \draw[->] (3.05,-3.5)--(4.5,-3.5);
  \node at (4.2,-5) {(b)};
\end{tikzpicture}
\caption{Two local changes of a cubic core. (a) A Whitehead move changes the pairing $(1,2)\mid(3,4)$ to $(1,3)\mid(2,4)$ through contraction; thick edges are the old and new middle edges. (b) A cubic vertex $v$ expands to a triangle, with $x_i$ joined to the original outside neighbor $r_i$. Deleting the dashed edge $x_2x_3$ leaves a subdivision of the old graph. Open circles mark ports to the omitted outside graph.}
\label{fig:moves}
\end{figure}
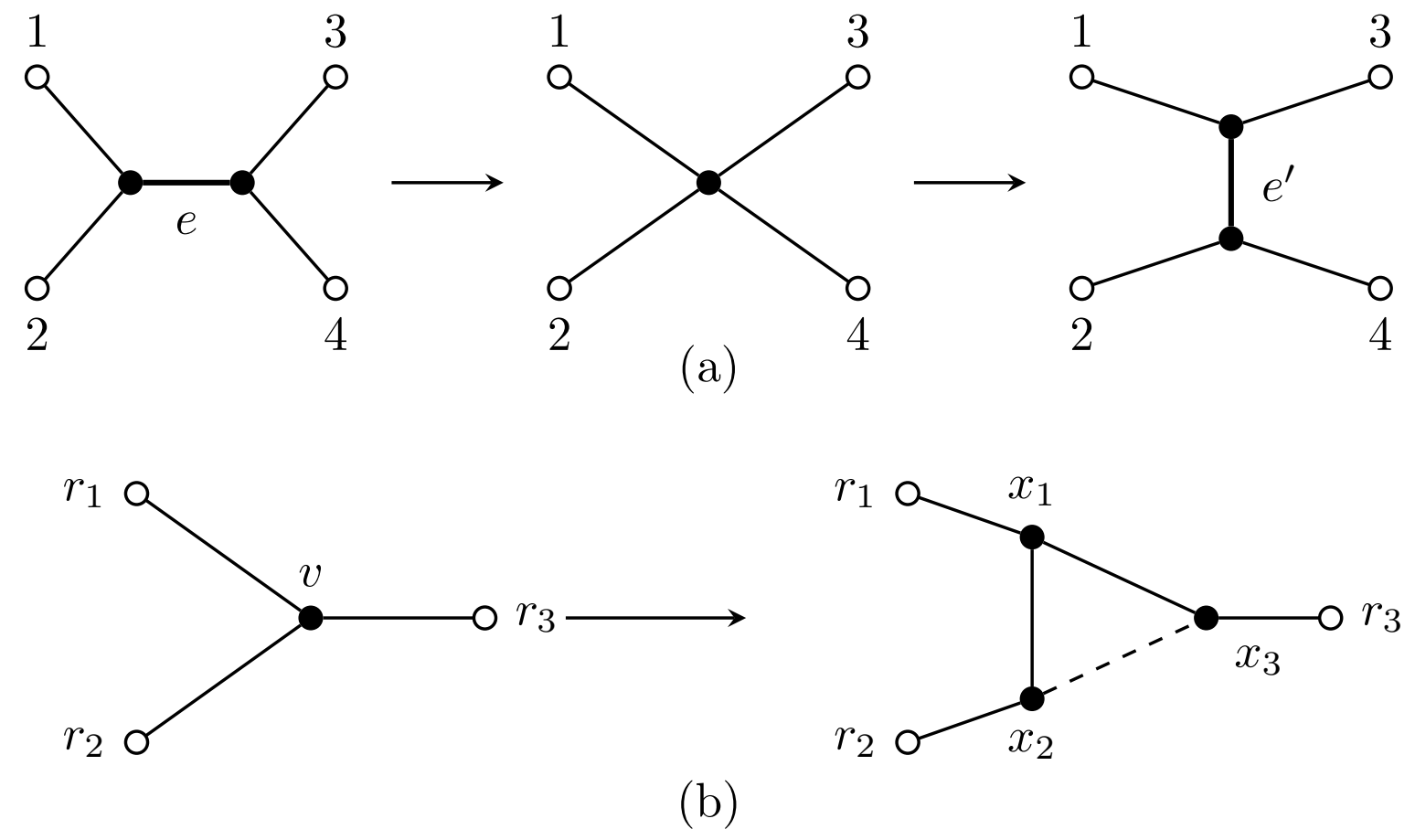

A Whitehead move contracts an edge between two cubic vertices and expands
the resulting four-valent vertex with another pairing of its four ports
(Figure~\ref{fig:moves}(a)).
If both resulting graphs are simple three-connected cubic graphs, then
\begin{equation}\label{eq:whiteheadall}
 \mathsf P(G)\quad\Longleftrightarrow\quad\mathsf P(G').
\end{equation}
To compare the data, gauge the middle edge to one. Its endpoint gauges
must then agree, so both gauge quotients identify with that of the same
contracted graph. It remains to determine each new simple-cycle trace
from the old ones. A cycle avoiding both new vertices is unchanged.
If it uses the new middle edge, removing that edge leaves one outside
path: its old ports either meet the same old vertex or are joined by the
old middle edge, giving an old simple cycle with the same trace. If it
visits only one new vertex, close its outside path with the old middle
edge. That edge has unit transport, so again the traces agree. The only
remaining case visits both new vertices without using the new middle edge. Label the old ports so that $1,2$ meet one old
vertex and $3,4$ meet the other, while the new pairing is $1,3$ and $2,4$.
The two outside arcs of this exceptional cycle are internally disjoint.
Since they join the two new vertices into one cycle, their endpoint pairs
are either $(1,2),(3,4)$ or $(1,4),(2,3)$.

In the first case, the old pairing closes the arcs into two disjoint
cycles joined by the old middle edge. The new cycle trace is the trace
of their joined dumbbell, which Lemma~\ref{lem:dumbbellall} determines.
In the second case, both arcs join the two old vertices. Together with
the old middle edge they form a theta graph. With that edge gauged to one,
orient both arcs from the first old vertex to the second and denote their
transports by $X,Y$. The old simple cycles give $\tau(X),\tau(Y)$ and
$\tau(XY^{-1})$, while the new cycle requires $\tau(XY)$, up to inversion.
The identity $\tau(XY)+\tau(XY^{-1})=2\tau(X)\tau(Y)$ supplies it.
The same argument applies to the reverse move. Thus equality of simple-cycle
data transfers in both directions under a bijection of gauge quotients,
proving \eqref{eq:whiteheadall}.

An open ear is a path with distinct old endpoints and new internal
vertices; a missing edge is an ear of length one. A biconnected subgraph
can be enlarged to the whole biconnected graph by ears: join two old
neighbors of an outside component, and finally add missing edges.
Each ear increases cycle rank by one.

The triangle expansion requires an additional observation about low-rank
path data. Their apparent sphere of solutions is removed by the old
stabilizer. The following lemma identifies that stabilizer without a
nondegeneracy assumption.

\begin{lemma}\label{lem:pathsubgroup}
Let $H$ be biconnected and let $a\ne b$. Normalize one simple $a$--$b$
path transport to one. If all such transports belong to a subgroup
$K\le SU(2)$, the entire connection is gauge equivalent to a $K$-connection.
If the old connection is central or commuting noncentral, the real path
span has dimension one or two, respectively. Adding an ear then gives
at most one new gauge class for each fixed old class and its new cycle data.
\end{lemma}
\begin{proof}
First order the vertices with $a$ least and $b$ greatest, so that every
other vertex has a smaller and a larger neighbor. Begin with an $a$--$b$
path. Each component outside the ordered set has two distinct neighbors
in it, since $H$ has no cut vertex. Insert a path through that component
between two such neighbors and repeat. Orient all edges upwards.
Increasing paths from $a$ to any vertex and from that vertex to $b$ now
exist. Choose such prefixes and suffixes with transports $\ell_x,r_x$.
For an increasing edge $xy$, every vertex of its chosen prefix is at most
$x$, whereas every vertex of its suffix from $y$ is at least $y>x$.
Thus their concatenation is simple. The same observation applies to the
prefix and suffix meeting at $y$. Both $\ell_xT_{xy}r_y$ and
$\ell_yr_y$ therefore belong to $K$, and their quotient is
$\ell_xT_{xy}\ell_y^{-1}$. Taking the vertex gauges to be $\ell_x$
puts every upward edge, and hence its inverse, in $K$.

If the old connection is central, its normalized endpoint paths lie in
$\{\pm1\}$, so their real span has dimension one. If it is commuting,
they lie in a common circle subgroup and span at most two dimensions.
Conversely, the normalized span contains one. A one-dimensional span
has unit section $\{\pm1\}$, while a two-dimensional span is generated
by one and a purely imaginary vector and has a circle subgroup as its
unit section. The first assertion then makes the old connection central
or commuting, respectively. In particular, a commuting noncentral
connection has path span exactly two.
In a compatible tree gauge the central and commuting types have effective
stabilizers $SO(3)$ and $SO(2)$, respectively. The new total ear transport
$q$ and each old simple endpoint path $P$ are oriented in the same direction.
Write $p_P$ for its normalized transport and $c_P$ for the prescribed
half-trace of the cycle formed by the ear and the reverse of $P$.
The ear equations are $\langle q,p_P\rangle=c_P$ and $|q|=1$. These equations fix its projection onto the real path span.
For a central old connection only the scalar part of $q$ is fixed;
its remaining imaginary vector lies on a sphere in $\R^3$, on which the
old $SO(3)$ stabilizer is transitive. In the commuting noncentral case,
the scalar part and the component along the common imaginary axis are
fixed. The remaining two components lie on a circle, and the old
$SO(2)$ stabilizer rotates that plane transitively. Radius zero gives a
single point in either case; an incompatible radius gives no extension.
Thus every nonempty solution set is one orbit of the actual old stabilizer.
\end{proof}

\begin{lemma}\label{lem:triangleall}
Expanding a vertex into a triangle preserves $\mathsf P$ among simple
three-connected cubic graphs.
\end{lemma}
\begin{proof}
Let $G_0$ have $\mathsf P$, and expand its vertex $v$ to $x_1,x_2,x_3$,
with outside neighbors $r_1,r_2,r_3$ (Figure~\ref{fig:moves}(b)).
Write $R=G_0-v$. It is connected (indeed biconnected), because deleting
$v$ and any further vertex leaves $G_0$ connected; the three outside
neighbors are distinct because $G_0$ is simple and cubic.
Deleting $x_2x_3$ leaves a subdivision $H$ of $G_0$. Align the two
connections on $H$, using $\mathsf P(G_0)$, and gauge $x_1x_2,x_1x_3$
to one. Only $q=T_{x_2x_3}$ may differ from $q'$. Write $P_i=T_{x_ir_i}$.
If $H$ is commuting, Lemma~\ref{lem:pathsubgroup} finishes the proof.

It remains to treat a noncommuting $H$. The outside graph $R$ can have
a larger stabilizer than $H$, so we separate the central and noncentral
cases. If $R$ is central, gauge its edges to
$\pm1$ and arrange $P_1=1,P_2=B,P_3=C$. The variables $B,C$ do not commute.
Choose a simple path in $R$ between each pair of outside neighbors.
The triangle-tree path $x_2x_1x_3$ has transport one. Consider, respectively,
\[
 x_2r_2\leadsto r_1x_1x_3,\qquad
 x_2x_1r_1\leadsto r_3x_3,\qquad
 x_2r_2\leadsto r_3x_3
\]
where each $\leadsto$ denotes the chosen path in $R$. These give three
further simple paths, with transports $B,C^{-1},BC^{-1}$
up to signs: all transports inside the central $R$ are $\pm1$.
Since $B,C$ do not commute, their imaginary parts are linearly independent.
The first three quaternions span one and these two imaginary directions;
the fourth has a nonzero component along their cross product. Thus these
four paths span $\R^4$. Their cycle traces fix all four inner products
of $q$ and $q'$, so $q=q'$.

If $R$ is noncentral, gauge its tree to one. A gauge preserving these
tree edges is constant on the connected $R$; preserving all its edges
requires this common value $h$ to centralize its chord tuple. Its
centralizer $Z$ is a circle subgroup when that tuple is commuting
noncentral, and $\{\pm1\}$ when it is noncommuting.
For $A\in SU(2)$ let $C_{SU(2)}(A)=\{g\in SU(2):gA=Ag\}$ be its
centralizer. Deleting $x_1x_3$ and contracting the other two triangle edges gives
outside transports $(P_1,P_2,qP_3)$ and $(P_1,P_2,q'P_3)$ on $G_0$.
They are gauge equivalent. Keeping the first two transports forces
\[
 h\in Z\cap C_{SU(2)}(P_1^{-1}P_2),\qquad
 k=P_1hP_1^{-1}=P_2hP_2^{-1}
\]
at the contracted vertex. The third transport satisfies
$kqP_3h^{-1}=q'P_3$. If the displayed intersection is central, then
$h=k=\pm1$, so this equation reduces to $q=q'$. Deleting $x_1x_2$ instead yields
$(P_1,q^{-1}P_2,P_3)$ and the same conclusion if
$Z\cap C_{SU(2)}(P_1^{-1}P_3)$ is central.
If both intersections are noncentral, then $Z$ must be a circle $T$.
The centralizer of any noncentral element of $T$ is $T$ itself. Hence
$P_1^{-1}P_2,P_1^{-1}P_3\in T$. Every tree-gauged edge of $R$
commutes with $T$ and therefore also lies in $T$. Apply the common
gauge $P_1^{-1}$ at the three triangle-tree vertices, keeping $R$ fixed.
The three outside transports become $1,P_1^{-1}P_2,P_1^{-1}P_3$,
while the triangle-tree edges remain one. All old transports now lie
in $T$, so $H$ is commuting, a contradiction.
\end{proof}

\subsection{Generating the cores and excluding the other graphs}
We now connect the local operations to the whole graph class. The generation
step initially inserts an edge between two subdivided edges. Sliding one
of them reduces this insertion to a triangle expansion, with Whitehead
moves providing the intervening comparisons.

Every simple three-connected cubic graph is generated from $K_4$ by
repeatedly subdividing two distinct edges and joining the new vertices
\cite[Section 4]{OkThomassen}. Denote this operation by $J(f,g)$.
Adjacent $f,g$ give a triangle expansion. Otherwise fix $g$ and slide
$f$ along a shortest edge chain until it first becomes adjacent to $g$.
For $f=uv,f'=vw$, with subdivision point $a$ on $f$ and point $b$ on $g$,
a Whitehead move on $av$ replaces $ua,vw$ by $uv,aw$, preserving $av,ab$
and the third edge at $v$. This gives exactly $J(f',g)$.
Each $J(f,g)$ is simple and cubic because $f,g$ are distinct edges of a
simple cubic graph. To check three-vertex-connectivity, consider a pair
of deleted vertices. If both are old, the surviving old graph is connected.
A new vertex can lose both old neighbors only when its subdivided edge
joins the two deleted vertices. The other new vertex still has an old
neighbor, since $f\ne g$, so the joining edge connects any stranded vertex
back to the old graph. If one old vertex and one new vertex are deleted,
the old graph minus that old vertex is biconnected and hence bridgeless.
Deleting the edge subdivided by the removed new vertex therefore leaves
the surviving old vertices connected; the other new vertex has at least
one surviving old neighbor. If both new vertices are deleted, the old
graph loses only $f,g$ and remains connected by three-edge-connectivity.
Deletion of fewer vertices also preserves connectivity. Thus all slides
satisfy \eqref{eq:whiteheadall}.

A shortest edge chain to $g$ can be stopped at its first edge adjacent
to $g$, so every joining operation along the chain uses distinct edges.
At that last edge, $J$ is a triangle expansion of the old core.
Consequently Whitehead equivalence transfers $\mathsf P$ from that
triangle expansion back along the chain to $J(f,g)$. Starting from
Lemma~\ref{lem:k4}, Lemma~\ref{lem:triangleall} and induction over the
quoted sequence of joining operations prove $\mathsf P$ for every core.

\begin{proof}[Proof of Theorem~\ref{thm:allconnections}]
Proposition~\ref{prop:equivalence} identifies the first two conditions.
\emph{Necessity.} Two cyclic blocks already fail: in a common unit tree
gauge choose one chord in each and assign $(i,i)$ and $(i,j)$, with all
other transports one. Simple cycles stay within one block and give the
same traces, whereas common conjugation cannot separate two identical elements.

For a vertex set $S$, let $\delta(S)$ be the set of edges with one endpoint
in $S$ and the other outside $S$. The remaining obstruction is such a cut
$\delta(S)$ of at least four edges,
whose two sides are connected and which every simple cycle crosses at
most twice. Orient the cut outwards and put $(1,i,j,k)$ and $(1,i,j,-k)$
on four cut edges, with all other transports one. Every cycle sees at
most two of these variables, so their common real Gram matrix gives
equal traces. Indeed, a cycle crossing the cut twice has half-trace
the real inner product of its two outward cut transports; all edges on
either side have unit transport. Choose a spanning tree on each connected
side and connect them by the cut edge assigned one in both connections.
This is a common unit spanning tree. Any gauge between the two
connections must therefore be constant on the whole graph, giving common
conjugation. A conjugation fixing $i,j$ also fixes $k=ij$, so it cannot
send $k$ to $-k$.

In a cyclic block $B$, a vertex of degree at least four gives this cut
with $S$ the singleton vertex, since $B-v$ is connected. Hence all degrees
are two or three. The all-degree-two case is a cycle. Otherwise suppress
degree-two chains, obtaining a loopless biconnected cubic multigraph $K$.
A loop would represent a cycle attached to the rest through its single
cubic endpoint, making that endpoint a cut vertex before suppression.
On each suppressed path, internal vertex gauges set all but its final
transport to one, leaving precisely its total path transport. Thus
suppression preserves gauge classes; simple cycles correspond under
suppression, including two-edge cycles in $K$. The same observation
allows a mirror obstruction constructed on $K$ to be lifted to $B$.

If $K$ has a two-edge cut $ac,bd$, write its sides as $A,B$, with
$a,b\in A$ and $c,d\in B$. Neither side is a singleton: a single vertex in a loopless cubic graph
has boundary three, whereas this cut has size two.
The endpoints on either side are distinct: if $a=b$, deleting $a$
separates the nonempty $A-a$ from $B$, contrary to biconnectivity; the
other side is analogous. Every component of $A-a$ must contain $b$,
since otherwise it would be separated from $B$ by deleting $a$.
Thus $A-a$ is connected, and similarly $B-c$ is connected.
The set $S=\{a,c\}$ is connected by $ac$ and has four boundary edges;
its complement is connected through $bd$. A simple cycle using $ac$
must also use $bd$, since those are the only edges of the original cut.
It then uses one additional incident edge at each of $a,c$, crossing
$\delta(S)$ exactly twice. A cycle using neither cut edge lies in one
side and visits at most one of $a,c$, so crosses $\delta(S)$ at most twice. The preceding mirror obstruction applies and survives
subdivision back to $B$.

If $K$ has two vertices, it is a triple edge and $B$ is a theta graph.
With more vertices, parallel edges give a two-edge cut at their endpoints,
so $K$ is simple. A two-vertex cut $\{u,v\}$ in a biconnected cubic graph
also gives a two-edge cut. Indeed, each remaining component meets both
vertices. If $uv$ is present, the four outside ports are distributed among at
least two components. Each component uses at least one port at each
of $u,v$, so there are exactly two components, both with boundary two.
Otherwise the six ports give two or three components. For three
components each uses one port at each cut vertex and has boundary two.
For two components, write $(r,s)$ for the numbers of edges from one
component $C$ to $u,v$. Both entries lie in $\{1,2\}$, and the other
component has counts $(3-r,3-s)$. If $(r,s)=(1,1)$, take $C$;
if $(r,s)=(2,2)$, take the other component. If $(r,s)=(2,1)$,
take $C\cup\{u\}$: its two outgoing edges are the one from $C$ to $v$
and the one from $u$ to the other component. For $(r,s)=(1,2)$,
take $C\cup\{v\}$ instead. Each case gives a two-edge cut, already
excluded. Hence the remaining core is three-vertex-connected.

\emph{Sufficiency.} Trees have one gauge class because there are no
chords after tree gauge. On a cycle the single remaining quaternion is
determined up to conjugacy by its scalar part, hence by its holonomy trace.
For a theta graph, normalize one terminal path to one and call the other
two transports $A=(a,u)$ and $B=(b,v)$. The three cycle traces are
$a,b$ and $ab+u\cdot v$; unit norms give $|u|^2=1-a^2$ and
$|v|^2=1-b^2$. Thus the imaginary Gram matrix is known. As its span
has dimension at most two, any Gram isometry can be made proper on the
orthogonal complement, giving one common $SO(3)$ orbit even at degenerate data. The cubic-core argument
above treats all other listed cores. Subdivision and attached trees add
only gauges that can be eliminated. This proves both directions.
\end{proof}

\begin{corollary}\label{cor:allcontinuous}
On every graph in Theorem~\ref{thm:allconnections}, $W_G$ is a homeomorphism
onto its image. Every continuous realizable data path therefore has a
unique continuous gauge-class lift, including at degenerate connections.
\end{corollary}
\begin{proof}
The domain is compact, the target is Hausdorff, and $W_G$ is continuous
and injective. Its inverse on its image is continuous.
\end{proof}

\section{Every fiber and its degenerations}\label{sec:everyfiber}
\subsection{Local choices and relative block rotations}
For arbitrary graphs, the data can leave finite local choices and
continuous relative block rotations. We first construct the local choices,
then describe how they assemble and meet at exceptional data. An ear adds
one quaternion subject to linear observations; the only continuous freedom
within a block is a stabilizer orbit. Continuous ambiguity in the final
fiber comes from placing different blocks relative to one another.

For the algorithmic statements, a real-algebraic number means a real root
of a nonzero polynomial over $\mathbb Q$, given exactly by such a polynomial
and a rational interval isolating the root. Real-algebraic input has such
coordinates. A subset of $\R^d$ is semialgebraic if it is a finite Boolean
combination of polynomial equalities and inequalities. An effective
semialgebraic description supplies these polynomials and their
real-algebraic coefficients explicitly.

On $(S^2)^a$, the group $SO(3)$ acts diagonally, rotating all $a$ vectors
by the same rotation. The quotient below uses this action.

\begin{theorem}\label{thm:everyfiber}
For a biconnected block $B$ of cycle rank $\gamma$, every realizable simple-cycle
datum has at most $2^{\max(\gamma-2,0)}$ gauge classes. All classes are
constructible by an ear recursion. Fix realizable data on a connected graph. Let $d_B$ be
these local fiber sizes, let $a$ count commuting noncentral blocks and
let $b$ count noncommuting blocks. The complete fiber has
$\prod_Bd_B$ connected components, each homeomorphic to
\[
 \begin{cases}
 (S^2)^a\times SO(3)^{b-1},&b\ge1,\\
 (S^2)^a/SO(3),&b=0.
 \end{cases}
\]
The types $a,b$ are constant across all local choices in a fixed fiber.
For real-algebraic input, the construction is a finite exact algorithm.
\end{theorem}
\begin{proof}
If $B$ is acyclic, tree gauge gives its unique class, with $\gamma=0$.
Otherwise start an open-ear decomposition of $B$ with a cycle of half-trace $c$,
whose unique class has representative $c+\sqrt{1-c^2}\,i$.
For each old candidate normalize one simple endpoint-path transport to
one. A new ear has one total variable $q\in S^3$; all new cycles are
that ear joined to an old simple endpoint path. For each old simple endpoint path $P$, orient it and the ear in the same
direction. The row of the real matrix $M$ indexed by $P$ is the four
coordinates of its normalized transport, and $y_P$ is the prescribed
half-trace of the new cycle obtained by returning along $P$.
The vector $y$ collects these observed half-traces. The complete extension
equations are
\begin{equation}\label{eq:allearequations}
 Mq=y,\qquad |q|=1,
\end{equation}
Check consistency of every row, take the minimum-norm solution $q_*$,
and intersect $q_*+\ker M$ with $S^3$. Since $q_*\perp\ker M$,
a solution $q=q_*+u$ has $u\in\ker M$ and
$|u|^2=1-|q_*|^2$. Thus the solution set is a sphere in the
$(4-\operatorname{rank}M)$-dimensional kernel, possibly empty or a singleton.
When $M$ has real rank one or two, Lemma~\ref{lem:pathsubgroup} applies: the unit solutions
form at most one orbit of the actual old stabilizer. When $M$ has real rank three or
four, the old connection is noncommuting and has trivial effective
stabilizer. There are at most two solutions in rank three, and one in
rank four. Negative squared radius gives no solution; in rank four the
unique solution must have unit norm.

Every actual extension restricts to a listed old class, can be aligned
with its representative, and then satisfies exactly these equations.
Conversely, put each retained $q$ on the last edge of the new ear and
put one on its other edges. This realizes every retained solution on
the original graph. Extensions of distinct old classes stay inequivalent
by restriction. The initial cycle has cycle rank one, and each added ear increases
the cycle rank by one. The first ear follows a commuting cycle and does
not double; each of the remaining $\max(\gamma-2,0)$ ears increases
the number of classes by at most a factor of two. This proves the bound.
If a final class is commuting, its restrictions are commuting at every
stage. Starting from the unique initial cycle class, induction then gives
a unique class at each stage: the old representative is commuting, so
Lemma~\ref{lem:pathsubgroup} leaves at most one extension class, and the
given final class guarantees existence. Thus the entire local fiber is a singleton. Thus a local
fiber is either one central class, one commuting noncentral class, or
a finite set of noncommuting classes.

Take the union of block trees and bridges as a global tree. For each
combination of local representatives $q_B$, every global class is obtained
by rotating the block tuples independently, modulo a common rotation.
The resulting component is
\begin{equation}\label{eq:allblockquotient}
 SO(3)\backslash\prod_B SO(3)/H_B,
\end{equation}
Here the product is over the cyclic blocks, and
$H_B=\{R\in SO(3):R\cdot q_B=q_B\}$ is the effective stabilizer of the
chosen block tuple. It is $SO(3)$, $SO(2)$ or $\{1\}$ in the three
respective types. The left quotient removes a common rotation of all blocks.
If there is a noncommuting block, its frame removes the common rotation;
otherwise only the axis quotient remains. Each quotient is connected.
Different local choices lie in disjoint compact subsets under restriction,
so these are exactly the connected components. All finite steps, equality
tests and square-root choices are real-algebraic operations.
\end{proof}

\subsection{Branch merging and stabilizers}
The preceding construction describes each fixed fiber. To compare nearby
fibers, we need coordinates that survive the loss of a frame. For a
tree-gauge tuple $q=(q_j)_{j=1}^{\beta}$, write $q_j=a_j+v_j$ with
$a_j\in\R$ and $v_j\in\R^3$, and use the complete invariants
\begin{equation}\label{eq:completeinvariants}
 I(q)=\bigl((a_j),(v_i\cdot v_j),(\det(v_i,v_j,v_k))\bigr).
\end{equation}
The coordinates include every scalar part, every pairwise inner product
and every triple determinant. They distinguish all gauge classes by the
Gram and orientation argument.
The stabilizer in $SU(2)$ is $SU(2)$, $U(1)$ or $\{\pm1\}$ when the
global imaginary span has dimension zero, one or at least two; here
$U(1)=\{z\in\C:|z|=1\}$ denotes the circle group.
Two limiting branches merge exactly when their complete global invariants
agree. This test includes cross-block inner products and oriented volumes;
equality of the separate local classes alone is insufficient.

\begin{proposition}\label{prop:degenerationatlas}
For each fixed graph there is a finite semialgebraic construction of the
whole realizable data image, all fibers, their branch incidences and their
stabilizer strata. For an effective semialgebraic data path, its restricted
total space and endpoint reachability are also finitely constructible.
\end{proposition}
\begin{proof}
Keep the unit quaternion variables in \eqref{eq:allearequations} and
project to the complete invariants \eqref{eq:completeinvariants}.
The ear recursion gives all finite block fibers. Effective semialgebraic
decomposition \cite[Sections 2 and 3.5]{Basu} partitions the data image
so that these finite fibers are disjoint continuous semialgebraic graphs.
Refine the partition to fix imaginary rank and nonzero frame pivots.
Set $G_{ij}=v_i\cdot v_j$ and $T_{ijk}=\det(v_i,v_j,v_k)$.
For a pair of pivot indices $p,q$ with independent vectors, a representative is given,
in rank at least two, by
\[
 v_l=\left(\frac{G_{pl}}{\sqrt{G_{pp}}},
 \frac{G_{pp}G_{ql}-G_{pq}G_{pl}}{\sqrt{G_{pp}D_{pq}}},
 \frac{T_{pql}}{\sqrt{D_{pq}}}\right),\qquad
 D_{pq}=G_{pp}G_{qq}-G_{pq}^2>0.
\]
In rank one choose $p$ with $G_{pp}>0$ and use
$v_l=G_{pl}e_1/\sqrt{G_{pp}}$, where $e_1=(1,0,0)$; in rank zero use zero.
These formulas preserve all Gram entries and orientations.

For each combination of block graphs, form the joint representative graph
$\Gamma=\{(z,q(z))\}$ over a common data cell: $z$ is the Wilson data vector
and $q(z)$ is the tuple of the chosen block representatives at those data. Take its closure in the
compact product of the data image and $(S^3)^\beta$, multiply by the
independent block rotation groups, and map to $(z,I(q))$.
There are finitely many such compact sources. Their union maps onto the
entire graph of $W_G$, since every point is represented on its own cell;
all closure points remain valid by continuity. Identify precisely those
source points with the same $z$ and the same complete global invariants.
The compact-to-Hausdorff quotient is exactly the full fiber family.
Closure membership, the identification relation and stabilizer rank are
semialgebraic conditions, so their boundaries and intersections are
finitely computable. Joint closure keeps the same approaching data in
all blocks and prevents incompatible marginal limits from being combined.

For a specified effective semialgebraic path, do this construction after
restricting to that path. Its time axis has a finite subdivision into open
intervals with disjoint local branches and the endpoint times. The displayed
representatives are bounded semialgebraic functions on each interval and
therefore have endpoint limits. Each branch combination supplies a closed
interval times the internal quotient \eqref{eq:allblockquotient}.
At an endpoint its stabilizers can increase; map the endpoint quotient
to the actual limiting tuples. This map covers the entire component of
that local endpoint combination. Include every vertex-time fiber and
identify endpoint images by their complete global invariants. The same
compact quotient argument gives the restricted total space.
Endpoint reachability is finite graph reachability among these local
branches, with the prescribed endpoint classes. Within one interval,
any two permitted endpoint configurations can be joined by choosing paths
in the connected group $SO(3)$ for the block rotations.
\end{proof}

The path restriction is essential: the closure of a full data cell contains
approaches that a prescribed path may never follow. Arbitrary continuous
paths are treated separately below.

\section{Continuous data and prescribed endpoints}\label{sec:continuousdata}
On the rigid graphs, Corollary~\ref{cor:allcontinuous} gives a unique
continuous lift of every realizable data path. On other graphs, local
branches may be forced to switch without meeting. The next example shows
that pointwise realizability alone does not resolve this obstruction.

\begin{example}[Realizable data with no continuous lift]\label{ex:nolift}
Write $K_{p,q}$ for the complete bipartite graph with parts of sizes $p,q$.
Take $K_{2,4}$ with terminals $s,t$ and internal vertices $x_1,\ldots,x_4$,
and add $x_1x_2$. In right-ordered transport coordinates put one on $sx_j$,
and put $1,i,j,\sigma k$ on $x_jt$, where $\sigma\in\{1,-1\}$.
Put $a+gj$ on $x_1x_2$, with $a=\sqrt{1-g^2}$ and $0\le g<1$.
The six old simple-cycle half-traces are zero. The six new ones are obtained by following $x_1x_2$ and returning
along the old simple $x_1$--$x_2$ paths. In the order
\[
 x_1sx_2,\quad x_1tx_2,\quad x_1sx_3tx_2,\quad
 x_1sx_4tx_2,\quad x_1tx_3sx_2,\quad x_1tx_4sx_2,
\]
the path transports are $1,-i,k,-\sigma j,-j,-\sigma k$.
The new cycle half-traces are therefore
\[
 (a,0,0,-\sigma g,-g,0).
\]
The old imaginary triple volume, measured in the ordered Hamilton basis
$(i,j,k)$, is $\sigma$; the two old classes never merge. For $0<t\le1/2$, put $f(t)=t\sin(1/t)$, $g(t)=|f(t)|$ and
set $\sigma g=f(t)$ in the displayed data, extending by $a(0)=1,g(0)=0$.
This defines continuous realizable data. Wherever $f\ne0$, it forces
$\sigma=\operatorname{sign}f$. A continuous lift would have a continuous
old triple volume taking these alternating values in $\{1,-1\}$, which
is impossible. Taking $f(t)=t$ instead makes only the positive-volume
endpoint reachable, although both classes belong to the endpoint fiber.
\end{example}

For a general continuous path, the number of degenerations can be infinite,
so the finite semialgebraic subdivision from the preceding section need
not apply. We instead select classes on finer and finer time grids.
The criterion asks whether the largest necessary jump tends to zero.
Give each block quotient the Euclidean metric of its invariant image
\eqref{eq:completeinvariants}.

For the graph application, $K$ is the compact invariant image of a block
quotient under \eqref{eq:completeinvariants}, and $Y$ is its realizable
Wilson image, a subspace of a finite-dimensional Euclidean space. The
induced observation map is $f:I(q)\mapsto W_B(q)$, where $B$ is the block.

\begin{lemma}\label{lem:gridlift}
Let $K\subset\R^N$ be compact, let $Y$ be a Hausdorff space, and let
$f:K\to Y$ be continuous with finite fibers, and let $z:[0,1]\to f(K)$ be continuous. Fix
$x_0\in f^{-1}(z(0))$. On the grid $j2^{-n}$ choose a point in each
fiber, starting at $x_0$, and let $b_n$ be the least possible largest
distance between consecutive choices. A continuous lift starting at
$x_0$ exists if and only if $\liminf b_n=0$, equivalently $b_n\to0$.
\end{lemma}
\begin{proof}
A continuous lift is uniformly continuous, proving necessity and the
full limit. For sufficiency choose grids of mesh tending to zero with
maximal jumps $e_l\to0$, and interpolate their chosen points linearly
in the ambient space.
Fix $\eta>0$ and a time $t_*$. Surround the finite fiber at $t_*$ by
disjoint open balls of radius less than $\eta/4$, with positive gaps
between their closures. Compactness of $K$ and continuity imply that
all fibers at nearby times lie in their union: otherwise compactness
gives a convergent subsequence of outside points, and continuity and the
Hausdorff property of $Y$ place its limit in that fiber, a contradiction.
Choose a smaller time neighborhood with closure in this neighborhood.
Finitely many smaller neighborhoods cover $[0,1]$. For sufficiently fine
grids their relevant adjacent nodes stay in the corresponding larger
neighborhoods, and $e_l$ is smaller than every ball gap. Consecutive nodes
therefore cannot change balls. Convexity keeps each interpolating segment
in the same ball. A Lebesgue number for the finite smaller cover now gives
equicontinuity, with oscillation less than $\eta$ on sufficiently short
time intervals. The interpolants are uniformly bounded, so
Arzel\`a--Ascoli gives a uniformly convergent subsequence.
Write $h$ for this uniform limit. For each fixed $t$, choose an endpoint
$s_l$ of the grid interval containing $t$ and let $x_l$ be the chosen
fiber point there. Then $s_l\to t$ and the distance from $x_l$ to the
interpolant at $t$ is at most $e_l$, so $x_l\to h(t)$. Closedness of $K$
gives $h(t)\in K$. Moreover, $f(x_l)=z(s_l)$ converges both to
$f(h(t))$ and to $z(t)$; uniqueness of limits in the Hausdorff space $Y$
gives $f(h(t))=z(t)$. All interpolants start at $x_0$, so $h(0)=x_0$.
\end{proof}
The same criterion holds with a prescribed final point
$x_1\in f^{-1}(z(1))$: fix both the first and last choices in every
grid chain when defining $b_n$. Necessity and the compactness proof
above then preserve both endpoints.

\begin{proposition}\label{prop:globallift}
A continuous realizable data path on a connected graph has a lift from
a prescribed global class if and only if its restriction to every cyclic
block has a lift from the corresponding local class. Thus the block
grid criteria of Lemma~\ref{lem:gridlift} give a necessary and sufficient
countable limit criterion. Effective semialgebraic paths have the finite
criterion of Proposition~\ref{prop:degenerationatlas}.
\end{proposition}
\begin{proof}
Restriction gives necessity. For sufficiency, lift each local class path
to a continuous tree-gauge tuple. We use the continuous curve-lifting
theorem for a compact Lie group acting orthogonally
\cite[Section 11.1.1]{RainerLifts}, applied to $SO(3)$ acting on the
imaginary vectors, with scalar coordinates fixed. The lift preserves
unit norms since they are invariants. Align its initial tuple with the
prescribed tuple by a fixed block rotation. Assign all block tuples to
their chords in one global tree gauge. This is a continuous global
connection with the required initial class and every required cycle trace.
There are finitely many blocks and each local fiber is finite, so applying
Lemma~\ref{lem:gridlift} to each block proves the criterion. If both global endpoints are prescribed and each local lift has the
corresponding local endpoints, first align the initial block tuples as
above. For each block choose a rotation taking its final tuple to the
prescribed final block tuple. Since $SO(3)$ is path connected, choose
a continuous path of rotations from the identity to that rotation and
apply it to the whole block tuple along the path. This keeps its initial
tuple fixed and preserves every block Wilson datum while aligning the
final tuple. The same global tree-gauge assembly then realizes both
prescribed global endpoints.
\end{proof}

\section{Generic rigidity and generic fibers}\label{sec:generic}
We now ask for uniqueness at a generic connection, while continuing to
compare it with every possible competitor. Rigid minors suffice in this
setting because the path frames needed to extend their reconstruction
fail only on a proper analytic set. The all-connection theorem requires
the stronger cubic structure precisely because it must also handle those
failures.

Generic global uniqueness
requires each true connection in such a generic set to have a singleton \emph{entire}
fiber, including comparison connections outside the generic set.
Let $D$ be the doubled triangle, with two parallel edges between each pair
of its three vertices. Its cycle rank is four. We write $H\preccurlyeq G$ when $H$ is obtained from $G$ by deleting
vertices or edges and contracting edges. These ordinary minors retain
parallel edges; a $D$ minor in a simple graph need not be a topological
copy. Write $c(G)=|\mathcal C(G)|$.

\begin{theorem}\label{thm:classification}
For a biconnected graph, generic global uniqueness holds if and only if
\begin{equation}\label{eq:criterion}
 \beta\le2\quad\text{or}\quad K_4\preccurlyeq G
                    \quad\text{or}\quad D\preccurlyeq G,
\end{equation}
Otherwise the
entire generic fiber has exactly two points. Moreover,
$c(G)\ge\beta(\beta+1)/2$, and for $\beta\ge3$ uniqueness is equivalent
to strict inequality.
\end{theorem}

We use $\coprod$ for disjoint union; the sphere quotient in the next
formula is again by a common rotation of all sphere factors.

\begin{theorem}\label{thm:fibers}
For a connected graph let $a$ be the number of cyclic blocks of rank one,
$b$ the number of cyclic blocks of rank at least two, and $r$ the number
of blocks of rank at least three containing neither $K_4$ nor $D$ as a
minor. The entire generic fiber is homeomorphic to
\begin{equation}\label{eq:fibers}
 \begin{cases}
 \displaystyle\coprod_{\nu=1}^{2^r}(S^2)^a\times SO(3)^{b-1},&b\ge1,\\[2mm]
 (S^2)^a/SO(3),&b=0.
 \end{cases}
\end{equation}
An empty product is a point. Generic uniqueness thus holds exactly for
trees and graphs with a single cyclic block satisfying \eqref{eq:criterion}.
For disconnected graphs, fibers are products over connected components.
\end{theorem}

The exceptional fibers are constructed in Section~\ref{sec:everyfiber}.
The full generic fiber in \eqref{eq:fibers} already includes singular
relative-axis configurations, even though the true connection is generic.

\subsection{Path spans and rigid minors}\label{sec:anchors}
\begin{lemma}\label{lem:pathspan}
In a biconnected graph $H$ of rank $\gamma$, between any distinct vertices
$s,t$ one can select $\min(\gamma+1,4)$ simple paths using only the graph
whose transports are generically linearly independent over $\R$.
\end{lemma}
\begin{proof}
Use the bipolar order constructed in the proof of
Lemma~\ref{lem:pathsubgroup}, with $s$ least and $t$ greatest, and orient
every edge upwards. Every edge belongs to a directed $s$--$t$ path,
and all these paths are simple.

For each directed path $P$, let $\chi_P\in\R^{E(H)}$ have coordinate one
on its edges and zero on all other edges. If an edge functional
$\omega$ has common value $c$ on all paths, choose a directed prefix and
suffix at each vertex with sums $A_v,B_v$. For each upward edge $u\to v$,
\[
 A_v+B_v=c=A_u+\omega_{uv}+B_v,
 \quad\text{so}\quad\omega_{uv}=A_v-A_u.
\]
Conversely potential differences have constant path sum. Their dimension
is $|V(H)|-1$, so the path differences span a space of dimension $\gamma$.
Fix one reference path $P_0$. The vectors $\chi_P-\chi_{P_0}$
span this $\gamma$-dimensional space. For $\gamma\ge3$, choose paths
$P_1,P_2,P_3$ whose differences are independent, and choose real edge
coefficients $\xi_e^{(j)}$ such that
\[
 \sum_e(\chi_{P_i}-\chi_{P_0})_e\xi_e^{(j)}=\delta_{ij}
 \qquad(1\le i,j\le3).
\]
Set $\xi_e=\xi_e^{(1)}i+\xi_e^{(2)}j+\xi_e^{(3)}k$. The transports for
$U_e(\epsilon)=\exp(\epsilon\xi_e)$ have expansions
$1+\epsilon\sum_{e\in P}\xi_e+O(\epsilon^2)$. Subtracting the reference
row from three others gives determinant $\pm\epsilon^3+O(\epsilon^4)$.
For ranks one and two use the analogous smaller minor. A nonzero analytic
function on the connected edge manifold has a zero set with empty interior
and Haar measure zero. Endpoint gauges act orthogonally on all path
transports, preserving independence.
\end{proof}

\begin{lemma}\label{lem:ears}
Generic entire-fiber rigidity of a biconnected subgraph of rank at least
three persists under ears. Four graph-selected cycle traces per ear suffice.
\end{lemma}
\begin{proof}
Use the four old endpoint paths from Lemma~\ref{lem:pathspan}; each joins
the ear to form a simple cycle. Gauge any competitor to agree on the
old graph. Internal ear gauges leave only its total unit transport $Q$
unknown, and the four new traces give its inner products with the old
path transports. Independence determines $Q$ among all competitors.
Aligning total transports aligns the ear by internal gauges. Finitely
many nonzero analytic conditions have a common generic intersection.
\end{proof}

On $D$ orient the three sides as $0\to1$, $1\to2$, $2\to0$.
In a tree gauge the reference transports are $I,I,D_0$, and the
alternative transports are $A,B,D_0C$, respectively. The four unit
quaternions $D_0,A,B,C$ are arbitrary. Ten actual cycles give
\begin{equation}\label{eq:ddata}
 \begin{gathered}
 \tau(A),\tau(B),\tau(C),\tau(D_0),\tau(D_0A),\tau(D_0B),\tau(D_0C),\\
 \tau(D_0BA),\tau(D_0CA),\tau(D_0CB).
 \end{gathered}
\end{equation}
Three cycles use a local parallel pair; the others traverse all three
sides with the indicated choices of nonreference edges. Subdivision
replaces edge transports by path transports without changing the argument.

\begin{lemma}\label{lem:double}
The observations \eqref{eq:ddata} have a singleton entire generic fiber.
\end{lemma}
\begin{proof}
Write $D_0=(d,w)$, $A=(a,u)$, $B=(b,v)$, $C=(c,z)$, and set
\begin{align*}
 F_1&=w\cdot(u\times v),\qquad F_2=w\cdot(v\times z),\\
 F_3&=(d^2-|w|^2)w\cdot(u\times z)
 -2d\bigl(|w|^2u\cdot z-(w\cdot u)(w\cdot z)\bigr).
\end{align*}
Assume $F_1F_2F_3\ne0$. Put $w$ on the positive third axis,
$d=\cos\phi$, $|w|=\sin\phi>0$ with $0<\phi<\pi$, and the transverse part of $u$ on the
positive first axis. The first seven observations determine the scalar
parts, third components and positive transverse lengths. Let the
remaining transverse angles of $B,C$ be $\theta_B,\theta_C$.
For $X$ preceding $Y$ in the list $A,B,C$, write $x_0,y_0$ for their
scalar parts, $r_X,r_Y$ for their transverse lengths, and $z_X,z_Y$ for
their third imaginary coordinates. Let $\theta_X,\theta_Y$ be their
transverse angles, with $\theta_A=0$. Multiplication gives
\begin{align*}
 \tau(D_0YX)={}&d(x_0y_0-z_Xz_Y)-|w|(y_0z_X+x_0z_Y)\\
 &-r_Xr_Y\{d\cos(\theta_Y-\theta_X)+|w|\sin(\theta_Y-\theta_X)\}.
\end{align*}
Thus the last observations determine
$\cos(\theta_B-\phi)$, $\cos(\theta_C-\phi)$ and
$\cos(\theta_C-\theta_B-\phi)$.
Every competitor enters the same gauge and must independently choose
$\theta_B$ or $2\phi-\theta_B$, and $\theta_C$ or $2\phi-\theta_C$.
Reflecting only the first requires an unchanged angle or
$\sin(\theta_C-2\phi)=0$; reflecting only the second requires an
unchanged angle or $\sin\theta_B=0$; reflecting both requires
$\sin(\theta_C-\theta_B)\sin\phi=0$. In this gauge the three factors are
\[
 F_1=|w|r_Ar_B\sin\theta_B,\qquad
 F_2=|w|r_Br_C\sin(\theta_C-\theta_B),\qquad
 F_3=|w|r_Ar_C\sin(\theta_C-2\phi).
\]
Their nonvanishing excludes the three genuine alternatives above.
If a reflection fixes its angle modulo $2\pi$, it produces the same
quaternion and hence no additional class.

Nonemptiness follows from $D_0=(3/5,0,0,4/5)$, $A=(0,1,0,0)$,
$B=(0,0,1,0)$ and $C=(0,5/13,12/13,0)$: the three factors are
$4/5,-4/13,-816/1625$. Their nonzero analytic product defines the
required generic set.
\end{proof}

\begin{lemma}\label{lem:minor}
Let $H$ be a loopless biconnected multigraph of rank at least three with
generic entire-fiber rigidity. Every biconnected simple graph containing
$H$ as an ordinary minor has generic entire-fiber rigidity.
\end{lemma}
\begin{proof}
Choose one cross edge per edge of $H$. Within each branch set take a
minimal tree joining its cross-edge endpoints, and let $F$ be their
union with the cross edges. Contraction preserves rank. The graph $F$ is biconnected. To see this, delete a vertex $x$ in
the branch tree representing $v\in V(H)$. The other branch trees and
the cross edges between them form a connected graph, since $H-v$ is
connected. Every nonempty component of the affected tree minus $x$
contains a cross-edge endpoint: otherwise that component could be pruned,
contradicting minimality of the tree joining all such endpoints.
Each component therefore attaches to the connected graph outside that
tree. If the affected tree is a singleton, no component remains and
the same outside graph proves connectedness.

Gauge the branch forest to identity. The remaining connections and gauges
are exactly those of $H$. Each simple cycle of $H$ lifts to a simple
cycle of $F$, using the unique path in each visited tree. Two-edge cycles
also lift, with length at least three in a simple original graph.
Thus the required observations recover the connection on $F$.
Conditional on forest transports, cross edges remain independent Haar
variables; the gauge coordinates also preserve openness and density.
The good minor set therefore pulls back to a generic set. Add ears from
$F$ to the whole graph and apply Lemma~\ref{lem:ears}.
\end{proof}

\subsection{Gram ambiguity and the classification}\label{sec:negative}
A two-terminal network is a graph with two designated distinct vertices.
Parallel composition identifies the corresponding terminals of otherwise
disjoint networks. Series composition identifies the second terminal of
one network with the first of another, leaving the two other terminals.
An ordinary prefix or suffix is a path attached at one terminal, sharing
no other vertex with the network.

We use the following consequence of Zaslavsky's frozen series-parallel
construction \cite[Example 6.1]{Zaslavsky}. A biconnected graph excluding
$K_4$ and $D$ as ordinary minors has a two-terminal construction from
paths by parallel composition and by adjoining ordinary prefixes or
suffixes. In the frozen construction, an edge may be subdivided or
doubled only while it is unfrozen; doubling an edge freezes every other
edge on its root-to-root path.
To obtain this representation from the frozen construction, consider the
first doubled edge of a root path. Every other edge on that path is then
frozen, so those edges remain ordinary prefixes and suffixes of the
parallel pair. All later subdivisions and doublings take place inside
its branches, where the same argument applies recursively. Thus each
series composition has at most one cyclic piece. This argument uses the
construction of the original graph; it makes no Whitney flip of the
connection or its observations.

\begin{lemma}\label{lem:gram}
Such a network of rank $\gamma$ has $N=\gamma+1$ terminal paths.
Unordered pairs correspond bijectively to its simple cycles, whose
half-traces are the inner products of the path transports
$Q_1,\ldots,Q_N\in S^3$. Every observation fiber has one gauge class if
the Gram rank is at most three and two if it is four.
\end{lemma}
\begin{proof}
Under parallel composition terminal-path counts add and ranks add with
an extra number of components minus one. Ordinary prefixes and suffixes
change neither count, proving $N=\gamma+1$ by induction. Two paths have
a common prefix and suffix and differ through one pair of branches,
recursively within one component if necessary. Their differing portions
form a unique cycle. Conversely every cycle lies within one parallel
component or traverses two, giving the bijection. Common transports
cancel in $\operatorname{Re}(Q_iQ_j^{-1})$.

Modulo internal gauges with the endpoint gauges fixed, the terminal
transports are arbitrary complete coordinates. This is immediate for a
path and is preserved by parallel products. For a prefix or suffix, its old terminal becomes an internal vertex
of the enlarged network. Gauges at that vertex and the other new internal
vertices make the attached path identity while fixing the two new
endpoint gauges. The induction hypothesis then applies to the middle
network, with exactly the prescribed terminal transports. Thus the tuple
is both arbitrary and complete modulo internal gauges.
The two remaining endpoint gauges act by $Q_i\mapsto aQ_i b^{-1}$,
with the same $a,b\in SU(2)$ for every $i$; these transformations give
all of $SO(4)$. Gram data determine
one $O(4)$ orbit. It splits into two $SO(4)$ orbits at rank four; below
rank four an orthogonal-complement reflection repairs orientation while
fixing the tuple. All these possibilities are realized by the arbitrary
terminal coordinates.
\end{proof}

\begin{proof}[Proof of Theorem~\ref{thm:classification}]
If the graph is acyclic, its gauge quotient is a point and $c=\beta=0$,
so all assertions hold. Assume henceforth that it contains a cycle.
The $K_4$ and $D$ minors give uniqueness by Lemmas~\ref{lem:k4},
\ref{lem:double} and \ref{lem:minor}. Otherwise the quoted graph
construction and Lemma~\ref{lem:gram} apply. For $\beta\le2$ the rank
is at most three, so every fiber is a singleton. For $\beta\ge3$ the
arbitrary terminal tuple generically spans four dimensions and gives
exactly two classes. A competitor with the same Gram data automatically
has the same rank, so exceptional competitors are included.

The affine-path proof of Lemma~\ref{lem:pathspan} gives at least
$\gamma+1$ endpoint paths on a biconnected graph of rank $\gamma$.
An ear with distinct endpoints has interior disjoint from the old
graph, so its union with each old simple endpoint path is a distinct
new simple cycle. Adding an ear raises the cycle rank from $\gamma$ to
$\gamma+1$ and adds at least $\gamma+1$ cycles, exactly the increment
\[
 \frac{(\gamma+1)(\gamma+2)}2-\frac{\gamma(\gamma+1)}2
 =\gamma+1
\]
in the proposed lower bound. Consequently, writing $c$ for the cycle
count and $\beta$ for the current rank, the difference
$c-\beta(\beta+1)/2$ is nondecreasing along an open ear decomposition.
A starting cycle has rank one and one simple cycle, hence difference zero. In the excluded class the path-pair bijection gives
equality. The anchors have $(\beta,c)=(3,7)$ and $(4,11)$; their minimal
minor models preserve rank and inject cycles, so their strict excess
persists through ears. Finally apply Proposition~\ref{prop:equivalence}.
\end{proof}

\begin{proof}[Proof of Theorem~\ref{thm:fibers}]
Take the union of block trees and bridges as a generating tree. Chord
variables fall into cyclic-block groups with one common conjugation
remaining. Since every simple cycle lies in a block, each choice of
compatible local gauge classes has global realizations
\[
 SO(3)\backslash\prod_B SO(3)/Z_B,
\]
where $Z_B\le SO(3)$ is the effective stabilizer of the block chord tuple,
as in \eqref{eq:allblockquotient}. A generic rank-one
block is noncentral with orbit $S^2$. A generic higher-rank block is
noncommuting with orbit $SO(3)$. Both compatible classes of a block with a two-point generic fiber are noncommuting, since commuting tuples would have terminal Gram
rank at most two. There are $2^r$ independent local choices.

If $b\ge1$, fix the first frame to remove the common action and obtain
$(S^2)^a\times SO(3)^{b-1}$. Otherwise the axis quotient remains. These
are homeomorphisms by compactness and the Hausdorff property. Distinct
local choices form disjoint compact sets, yielding the asserted disjoint
union, including singular axis configurations.
With two cyclic blocks, choose one noncentral based holonomy from
each and transport both to a common basepoint along the generating tree.
Write them as $(a,u)$ and $(b,v)$, with $u,v\ne0$. Independently rotating
the second block by $R\in SO(3)$ leaves every observed simple-cycle
trace fixed, whereas the half-trace of the based product is
$ab-u\cdot Rv$. This varies over the nondegenerate interval
$[ab-|u||v|,ab+|u||v|]$ and is invariant under global gauge.
It therefore distinguishes global gauge classes. A single
cyclic block contributes its local fiber; trees contribute none.
Independent component gauges give the disconnected product.
\end{proof}

\begin{example}
The graphs $K_4$ and $K_{2,4}$ both have rank three and have generic fibers
of one and two points. For generic data, two cycles meeting at a vertex have fiber
$(S^2)^2/SO(3)$, a closed interval, whereas two disjoint cycles have one
class for each compatible generic datum. Attaching trees preserves a fiber.
\end{example}

\section{Finite measurements and stable recovery}\label{sec:finite}
The preceding results use the entire weight family, or all simple-cycle
traces. We next reduce the measurement set. There are two distinct tasks:
preserving every fiber with finitely many full spectra, and reconstructing
a generic fiber stably from selected scalar traces. The measurement counts
below refer to these different observations.

A single fixed spectrum is generally insufficient: the gauge quotient
for the complete graph $K_5$ on five vertices has dimension $15$, whereas its ten eigenvalues with prescribed
sum supply at most nine real parameters. Apply the constant-rank theorem
to the smooth characteristic-coefficient map, whose trace coefficient is
fixed. This precludes generic uniqueness for that observation.

For $M\ge1$ settings $w^{(s)}\in(0,\infty)^E$, $1\le s\le M$, put
$\boldsymbol w=(w^{(1)},\ldots,w^{(M)})$ and define
its finite spectral observation by concatenation:
\[
 S_{\boldsymbol w}:X_G\longrightarrow\R^{2nM},\qquad
 S_{\boldsymbol w}(U)=
 \bigl(\lambda_j(w^{(s)},U)\bigr)_{1\le s\le M,\,1\le j\le2n}.
\]
Each block of $2n$ entries is ordered increasingly. The list of settings is
known and is kept in this order.

\begin{theorem}\label{thm:finite}
Let $G$ be connected with $\beta\ge2$, fix $0<a_0<b_0<\infty$, and put
$P=(a_0,b_0)^m$, $M=6\beta-5$. An open dense full-measure set of lists
$\boldsymbol w\in P^M$ preserves all fibers of $W_G$, simultaneously for
all connections. Concatenating the ordered spectra into $S_{\boldsymbol w}$,
there is a constant such that
\begin{equation}\label{eq:uniformwilson}
 \|W_G(U)-W_G(V)\|_\infty
 \le C_{G,\boldsymbol w}\|S_{\boldsymbol w}(U)-S_{\boldsymbol w}(V)\|_\infty
\end{equation}
for all $U,V$. For $\beta=1$, one may take $M=3$; trees require no data.
\end{theorem}
\begin{proof}
Tree gauge gives $X_G=(S^3)^\beta/SO(3)$. Scalar parts, imaginary inner
products and triple products separate all orbits, including rank-deficient
ones by orientation correction on the orthogonal complement. They realize
the quotient as a compact semialgebraic set of dimension $d=3\beta-3$:
noncommuting tuples have trivial stabilizer; commuting noncentral tuples
have only $\beta$ angle parameters; central tuples are finite. For rank
one, the dimension is one.

Collect the coefficients of $\Tr L_w(U)^\ell$, $1\le\ell\le n$, as
polynomials in $w$ into a vector $p(U)\in\R^N$, where $N$ is the total
number of collected coefficients. Its compact semialgebraic
image $Y$ has dimension at most $d$. Formula~\eqref{eq:isolation} gives
a fixed affine relation $W_G=Ap+c$, with $A$ a fixed linear map and $c$
a fixed vector. Conversely equal $W_G$ gives equal
all-weight spectra and hence equal coefficients. Thus $p$ and $W_G$
have exactly the same fibers. Moreover, at any fixed setting the
moments of orders $1,\ldots,n$ determine the reduced degree-$n$
characteristic polynomial by Newton identities: each reduced eigenvalue
is counted twice in $\Tr$. Thus equality of the evaluated coefficient
vectors is equivalent to equality of the full spectra at those settings.

If $Y$ is not a point, take the compact semialgebraic set
\[
 K=\overline{\{(p-p')/|p-p'|:p,p'\in Y,\ p\ne p'\}},\qquad \dim K\le2d.
\]
Each $q\in K$ is a nonzero coefficient vector for an $n$-tuple of
polynomials. Its common zero set in $P$ has dimension at most $m-1$.
Let $E_{\boldsymbol w}$ evaluate that tuple at the $M$ settings. The
incidence set $E_{\boldsymbol w}q=0$ has dimension at most
$2d+M(m-1)<Mm$ when $M=2d+1$. Its projection to $P^M$ is relatively
closed by compactness of $K$, and is lower dimensional. Its complement
is the desired open dense full-measure set. This is the dimension
mechanism of \cite{Sontag}, with singular semialgebraic strata included.

For a good list, $\sigma=\min_{q\in K}|E_{\boldsymbol w}q|>0$, and
\[
 |p(U)-p(V)|\le\sigma^{-1}|E_{\boldsymbol w}(p(U)-p(V))|.
\]
The right side contains spectral moments at fixed settings. If all spectra
lie in $[0,R]$, an ordered eigenvalue error $\epsilon$ gives moment error
at most $2n\ell\max(1,R)^{\ell-1}\epsilon$. Composing with $A$ proves
\eqref{eq:uniformwilson}. The secant closure is needed for this uniform
bound; separation of distinct pairs alone would give only injectivity.
For $\beta=1$ the same argument uses $d=1$ and $M=2d+1=3$;
for a tree $X_G$ is a point.
\end{proof}

\begin{theorem}\label{thm:sparse}
There is a set $\mathcal S_G$ of simple cycles determined by the graph
alone, of size zero for $\beta=0$, at most one for $\beta=1$, and at
most $4\beta-5$ for $\beta\ge2$, such that for generic $U$ and arbitrary $V$,
\[
 \tau_C(U)=\tau_C(V)\ (C\in\mathcal S_G)
 \quad\Longleftrightarrow\quad W_G(U)=W_G(V).
\]
\end{theorem}
\begin{proof}
On a block of rank $\gamma\ge3$ with neither a $K_4$ nor a $D$ minor, fix the first four
terminal paths. Observe their six pairwise inner products and the four
reference inner products of each remaining path, for $4\gamma-6$ cycles.
The reference Gram matrix $A$ is generically positive definite, and with
$v_i$ the inner-product column for path $i$, all missing entries are
\begin{equation}\label{eq:completion}
 \langle Q_i,Q_j\rangle=v_i^TA^{-1}v_j.
\end{equation}
Every competitor has the same positive definite $A$. Its four reference
transports therefore form a basis of $\R^4$, and its coordinates in this
basis are also $A^{-1}v_i$. Hence \eqref{eq:completion} holds for every
competitor, without a genericity assumption on that competitor, and
recovers the full terminal Gram matrix. Ranks one and two
use all one and three cycles, respectively.
For a block with a $K_4$ minor, lift the seven anchor cycles and add four per
ear, giving $7+4(\gamma-3)=4\gamma-5$. The doubled-triangle anchor uses
ten and gives $10+4(\gamma-4)=4\gamma-6$.
At each ear, the preceding anchor and ear data align any competitor on
the old graph; the four new traces then determine the total ear transport
by Lemma~\ref{lem:ears}. The path-independence conditions are imposed
only on the true connection. Lemma~\ref{lem:pathspan} selects the paths
from the graph before any transports are known.
With $a,b$ the block counts from Theorem~\ref{thm:fibers}, the union over
blocks has size at most $4\beta-3a-5b\le4\beta-5$ when
$\beta\ge2$. It preserves local fibers and all relative block rotations.
\end{proof}

For a selected length-$\ell$ cycle, let $F(t)=F_U(t)$ be the spectral-moment
polynomial from the proof of Proposition~\ref{prop:equivalence} for that
cycle. Put $N=\ell+1$ and choose
\begin{equation}\label{eq:cheb}
 \theta_j=\frac{(2j+1)\pi}{2N},\quad t_j=\frac{1+\cos\theta_j}{2},\quad
 a_j=\frac{(-1)^{N+1+j}}{N}\tan\frac{\theta_j}{2},\quad 0\le j<N.
\end{equation}
Let $T_N$ be the degree-$N$ Chebyshev polynomial, characterized by
$T_N(\cos\theta)=\cos(N\theta)$. Lagrange interpolation at the shifted roots $t_j$ of $T_N(2t-1)$ gives $F(0)=\sum_j a_jF(t_j)$.
Reversing the sum and using $\cot x\le1/x$ yields
\[
 \Lambda_N:=\sum_j|a_j|\le\frac4\pi\sum_{j=0}^{N-1}\frac1{2j+1}
 \le\frac4\pi\bigl(1+\tfrac12\log(2N-1)\bigr).
\]
For $\beta\ge2$, at most $(n+1)(4\beta-5)$ explicitly specified positive
complete spectra recover the selected cycle data. For $\beta=1$,
at most $n+1$ suffice; trees require none. This counts spectra, whereas
$4\beta-5$ counts scalar half-traces. If measured eigenvalue errors are
at most $\epsilon\le1$, the resulting trace error is at most
$D_{\mathcal S}\epsilon$. With $\Delta$ the maximum degree of $G$, define
\begin{equation}\label{eq:noisechain}
 D_{\mathcal S}=\frac n2\max_{C\in\mathcal S_G}
       \Lambda_{|C|+1}(2\Delta+1)^{|C|-1}.
\end{equation}
For an empty selection set $D_{\mathcal S}=0$.
True spectra at these settings lie in $[0,2\Delta]$.
The estimate retains the amplification from high-order moments.

\begin{remark}
On one biconnected block of rank at least two, any $C^1$ real scalar observation
map with finite generic fibers needs at least $3\beta-3$ coordinates,
by the constant-rank theorem. The bound $4\beta-5$ is not claimed optimal.
Its generic qualifier is necessary for this selection: the terminal tuples
$(1,1,1,1,i,i)$ and $(1,1,1,1,i,j)$ on $K_{2,6}$ agree on the fourteen
reference pairs and disagree on the omitted pair, whose values are one
and zero.
\end{remark}

\subsection{Stable recovery of complete fibers}\label{sec:stability}
We measure errors between the full fibers $\F_U\subset X_G$ defined in
Section~\ref{sec:model}. Fix a generating tree $T$ and chord representatives
$q,p\in(S^3)^\beta$.
Use the maximum chord quotient distance
\[
 d_{T,\infty}([q],[p])=\min_{R\in SO(3)}\max_j|q_j-R\cdot p_j|,
\]
Here $[q],[p]$ are the corresponding gauge classes and $R$ rotates their
imaginary parts. For nonempty compact subsets $A,B\subset X_G$, the
Hausdorff distance induced by $d_{T,\infty}$ is
\[
 d_H(A,B)=\max\left\{
 \sup_{a\in A}\inf_{b\in B}d_{T,\infty}(a,b),
 \sup_{b\in B}\inf_{a\in A}d_{T,\infty}(a,b)\right\}.
\]
 For a tree the quotient is a point, and the distance
is defined to be zero. The Euclidean product quotient metric differs by at most
a factor $\sqrt\beta$.

For either finite list of spectral settings in the next proposition,
write $S=S_{\boldsymbol w}$ for the fixed observation map. A
minimum-residual estimate minimizes the maximum coordinate discrepancy
between $S(x)$ and the measured spectral vector over $x\in X_G$.

\begin{proposition}\label{prop:regular}
If $G$ is generically rigid, fix settings from Theorem~\ref{thm:finite}
or the explicit selected-cycle settings. There is a generic regular locus
$\Omega_{\rm reg}\subset X_G$ such that for every compact $K\subset\Omega_{\rm reg}$,
\[
 d_{T,\infty}(x,y)\le C_K\|S(x)-S(y)\|_\infty,
 \qquad x\in K,\quad y\in X_G.
\]
A minimum-residual estimate with spectral noise $\epsilon$ therefore has
error at most $2C_K\epsilon$.
\end{proposition}
\begin{proof}
Let $H$ concatenate the nonleading coefficients of
$\det(zI-L_{w^{(s)}}(U))$ over the fixed settings $s$.
Use this analytic characteristic-coefficient map in place of ordered eigenvalues. The principal quotient stratum has analytic local
sections: rotate two nonparallel imaginary vectors into the positive
first coordinate axis and into the first--second coordinate plane with
positive second coordinate, respectively. For rank one use a noncentral
angle. The complement has Haar measure zero.
On every coordinate ball, maximal differential rank must equal the
quotient dimension. Otherwise the constant-rank theorem produces positive
dimensional local fibers on an open set, contradicting generic whole-fiber
uniqueness. Rank drops thus lie in proper analytic zero sets. Removing
them and intersecting with the generic whole-fiber uniqueness locus
gives $\Omega_{\rm reg}$, with full measure also after lifting through
gauge coordinates. On this locus both $H$ and $S$ have singleton whole
fibers: Theorem~\ref{thm:finite}, or Theorem~\ref{thm:sparse} and the
interpolation formulas above, identify those fibers with Wilson fibers.

Local inverse maps and a finite cover of $K$ give a uniform estimate for
close pairs. Pairs with $x\in K$ and $d_{T,\infty}(x,y)\ge\delta$ form
a compact set on which $|H(x)-H(y)|$ is bounded below, since the whole
fiber at $x$ is a singleton. The quotient diameter handles these far
pairs, including degenerate competitors. Characteristic coefficients
are Lipschitz in eigenvalues on a fixed spectral box. Residual minimization
adds two measurement errors. Eigenvalue crossings cause no difficulty.
\end{proof}

An effective real-algebraic procedure below is finite for rational or
real-algebraic input. Its output is a
parameterized set of classes, including every discrete branch and all
relative block rotations.

\begin{theorem}\label{thm:noisy}
For every connected graph and the cycle selection of
Theorem~\ref{thm:sparse}, there is a finite algorithm with the following
property. Given real-algebraic estimates $y_C$ and a nonnegative
real-algebraic error bound $\eta$ with $|y_C-\tau_C(U)|\le\eta$, it either reports an
uncertified outcome or produces a compact parameterized set
$\widehat\F\subset X_G$ and a computed bound $E_G$ satisfying
\begin{equation}\label{eq:hausdorff}
 d_H(\widehat\F,\F_U)\le E_G.
\end{equation}
For every true $U$ in a fixed generic set, all sufficiently small admissible
noises lead to success and $E_G=O_{U,G}(\eta)$. In combination with
\eqref{eq:noisechain}, this gives finite spectral reconstruction with
$O_{U,G}(\epsilon)$ fiber error.
\end{theorem}
The construction first certifies an inverse bound on each rigid anchor,
then propagates errors along the ears. In the Gram case it retains both
orientation branches and bounds their displacement. Independent block
rotations finally produce the whole approximate fiber. Appendix~\ref{app:stable}
gives the certification tests and the resulting error recursion.

The finite algebraic procedure may be expensive, and the bounds can grow
through a long ear chain. Numerical implementations need certified
solution and normalization errors if they claim rigorous certificates.
Disconnected graphs are handled componentwise with the maximum product
metric. None of these estimates asserts a uniform inverse constant
through every degenerate connection.

\section{Restricted edge controls}\label{sec:arbitrarycontrols}
Fix a controlled edge set $F\subset E$ and weights
$a=(a_e)_{e\in E\setminus F}\in(0,\infty)^{E\setminus F}$ on the other
edges. The control vector is $x=(x_e)_{e\in F}\in(0,\infty)^F$;
$L_{(x,a)}$ assigns weight $x_e$ to $e\in F$ and $a_e$ otherwise.
For these fixed data, we ask whether equality of spectra for every positive control
vector implies equality of all original simple-cycle Wilson traces.
This concerns Wilson observability; any gauge ambiguity inside a Wilson
fiber is still described by Theorem~\ref{thm:everyfiber}.

\begin{theorem}\label{thm:controltest}
For every finite labelled simple graph, control set $F$ and fixed positive
real-algebraic weights $a$, Wilson observability has a finite exact decision
procedure. Failure returns two actual connections with identical controlled
spectra and different Wilson data. Keeping $a$ symbolic gives an exact
semialgebraic condition for specified weights and for all positive weights.
\end{theorem}
\begin{proof}
Use a spanning forest and let $q,p$ be two tuples of chord quaternions.
Impose a unit-norm equation on each. For the entire graph, form
\[
 \chi_q(\lambda,x;a)=\det_{\C}(\lambda I-L_{(x,a)}(q)).
\]
Let $c_\nu(q,p;a)$ be every coefficient in $(\lambda,x)$ of
$\chi_q-\chi_p$, with $\nu$ indexing the monomials in those variables, and let
\[
 \Delta_W(q,p)=\sum_{C\in\mathcal C(G)}
                       (\tau_C(q)-\tau_C(p))^2.
\]
All are real polynomials in real quaternion coordinates. Observability
fails precisely when the following finite system has a real solution,
with an auxiliary scalar $r\in\R$:
\begin{equation}\label{eq:controlcollision}
 |q_j|^2=|p_j|^2=1,\qquad c_\nu(q,p;a)=0,\qquad
 r\Delta_W(q,p)-1=0.
\end{equation}
Indeed, an actual collision can be put in forest gauge, with the same
connection fixed for every control. Equality on the positive control
orthant is a polynomial identity, so all coefficients vanish. Unequal
Wilson data give $r=1/\Delta_W$. Conversely, any real solution gives
unit quaternion transports on the actual chords, equal characteristic
polynomials for all controls, and different original simple-cycle traces.

Effective real quantifier elimination and sampling \cite[Section 2]{Basu}
decide this system and return algebraic witnesses for algebraic input.
For all positive fixed weights, also make each $a_e$ unknown and impose
$a_es_e^2-1=0$ with auxiliary real variables $s_e$. A solution is exactly a failure at some positive weight.
Eliminating the quaternion and auxiliary variables instead returns the
parameter condition for individual real $a$.
Forests have $\Delta_W=0$ and thus an inconsistent last equation;
empty control sets, isolated vertices and disconnected graphs are included.
\end{proof}

For disconnected graphs one must retain the whole characteristic polynomial
in this test. Two disjoint unit-weight triangles with half-traces $(1,-1)$
and $(-1,1)$ have the same whole spectrum but different labelled Wilson
data. A triangle with half-trace $c$ has reduced polynomial
$\lambda^3-6\lambda^2+9\lambda-2+2c$; the two products coincide.
Thus separate observability of the components does not imply observability
from their unlabelled union of eigenvalues. The forest-gauge construction
uses that union directly.

\subsection{Rooted responses}\label{sec:rootresponses}
For attachments meeting the controlled core at one root, Schur elimination
makes the loss of information explicit. We fix the attachment weights and
vary specified core weights. Only a pair of scalar response polynomials
reaches the core, and different attachment Wilson data can have the same
response. We compute this pair, identify its kernel for triangle bouquets,
and then give core weights under which the controlled spectra recover it.

\subsubsection{Root elimination and gluing}
For a rooted graph $J$ with $m$ nonroot vertices, write $P_J(z)$ for its
monic degree-$(m+1)$ reduced characteristic polynomial, defined by
$\det_{\C}(zI-L_J)=P_J(z)^2$. Let $D_J(z)$ be the corresponding monic degree-$m$
polynomial of the root-deleted principal matrix, and put
\[
 C_J(z)=P_J(z)-zD_J(z).
\]
All principal deletions retain the original Laplacian degrees. A deleted
root edge of weight $a$ therefore still contributes $a$ to the diagonal
of its other endpoint. Away from zeros of $D_J$, the root Schur complement
is the scalar quaternion
\begin{equation}\label{eq:rootschur}
 z+\frac{C_J(z)}{D_J(z)},\qquad
 \frac{C_J(z)}{D_J(z)}=-\delta+O(z^{-1}),
\end{equation}
where $\delta=\sum_{e\ni\mathrm{root}}w_e$ is the weighted root degree in $J$. The Schur complement is scalar
because it is a one-by-one self-dual quaternionic matrix.
For real $z$ outside the deleted spectrum, block elimination gives
$\det_{\C}(zI-L_J)=D_J(z)^2$ times the square of this scalar
Schur complement. The monic normalizations, or equivalently their
asymptotics at infinity, select the identity
$P_J/D_J=z+C_J/D_J$ with no sign ambiguity.

Suppose a core $K$ and $J$ meet only at the root. Write $P_K,R_K$ for the
core reduced characteristic polynomial and root-deleted polynomial, again
with the original core degrees. Eliminating the nonroot variables of $J$
adds $C_J/D_J$ to the root entry of $zI-L_K$. Linearity of the reduced
determinant in that real diagonal entry gives
\begin{equation}\label{eq:rootgluing}
 P_{K\cup J}=D_JP_K+C_JR_K.
\end{equation}
These formulas first hold away from finitely many poles; multiplication
by the denominators makes them polynomial identities. Equal pairs
$(D_J,C_J)$ thus give equal spectra after gluing to any identical core,
for every choice of its allowed weights. This implication requires no
reconstruction hypothesis on the core.

\subsubsection{Triangle bouquets and the pole graph}
Let $J$ consist of $r$ triangles meeting only at the root. On petal $j$,
the two root edges have weights $\alpha_j,\beta_j>0$ and the remaining
edge has weight $\gamma_j>0$. Gauge the root edges to the identity and
write $Q_j$ for the internal transport and $\tau_j=\operatorname{Re}Q_j$.
Every vector $\tau\in[-1,1]^r$ is realizable: one may independently choose
\begin{equation}\label{eq:petalrealization}
 Q_j=\operatorname{diag}\bigl(\tau_j+i\sqrt{1-\tau_j^2},
                            \tau_j-i\sqrt{1-\tau_j^2}\bigr).
\end{equation}
Define
\begin{align*}
 H_j(z)&=(z-\alpha_j-\gamma_j)(z-\beta_j-\gamma_j)-\gamma_j^2,\qquad
 c_j=\alpha_j\beta_j\gamma_j,\\
 M_j(z)&=\alpha_j^2(z-\beta_j-\gamma_j)
                  +\beta_j^2(z-\alpha_j-\gamma_j).
\end{align*}
The deleted block and its inverse are
\[
 K_j=\begin{pmatrix}z-\alpha_j-\gamma_j&\gamma_jQ_j\\
 \gamma_jQ_j^*&z-\beta_j-\gamma_j\end{pmatrix},\qquad
 K_j^{-1}=\frac1{H_j}\begin{pmatrix}
 z-\beta_j-\gamma_j&-\gamma_jQ_j\\
 -\gamma_jQ_j^*&z-\alpha_j-\gamma_j
 \end{pmatrix}.
\]
Consequently the root correction is
$(\alpha_j,\beta_j)K_j^{-1}(\alpha_j,\beta_j)^T
=(M_j-2c_j\tau_j)/H_j$, a scalar. Thus
\begin{equation}\label{eq:petalresponse}
 D_J=\prod_j H_j,\qquad
 \frac{C_J}{D_J}=-\delta-\sum_j\frac{M_j}{H_j}+2R_\tau,
 \qquad R_\tau=\sum_j\frac{c_j\tau_j}{H_j},
\end{equation}
where $\delta=\sum_j(\alpha_j+\beta_j)$. The connection dependence is
exactly the rational function $R_\tau$.

Each $H_j$ has two distinct positive roots $l_j<u_j$: its discriminant is
$(\alpha_j-\beta_j)^2+4\gamma_j^2>0$, and both its root sum and product are
positive. Form a multigraph $\mathcal P$ whose vertices are the distinct
pole values and whose $j$th edge joins $l_j$ to $u_j$. Equal denominators
give parallel edges, which must be retained. Orient each edge from $l_j$
to $u_j$. Let $B$ be the signed incidence matrix: column $j$ is $-1$
at $l_j$, $+1$ at $u_j$, and zero elsewhere. Set
\[
 s_j=\frac{c_j}{u_j-l_j}>0,\qquad
 \mathcal L=B\operatorname{diag}(s_1,\ldots,s_r).
\]
A forest here excludes cycles of length two made of parallel edges.

\begin{proposition}[Pole-graph criterion]\label{prop:poleforest}
For fixed positive weights on a rooted triangle bouquet, the response
pair $(D_J,C_J)$ separates all realizable simple-cycle Wilson data if and
only if $\mathcal P$ is a forest. For every $\tau\in[-1,1]^r$, its entire
Wilson-coordinate response fiber is
\begin{equation}\label{eq:rootcubefiber}
 (\tau+\ker\mathcal L)\cap[-1,1]^r.
\end{equation}
If $\tau$ is interior to the cube, this fiber has dimension
$r-|V(\mathcal P)|+\kappa(\mathcal P)$, where $\kappa$ counts connected
components. The assertion concerns Wilson coordinates, not gauge classes.
\end{proposition}
\begin{proof}
Partial fractions give
\[
 \frac{c_j}{H_j(z)}=s_j\left(\frac1{z-u_j}-\frac1{z-l_j}\right).
\]
The vector of residues of $R_\tau$ is therefore $\mathcal L\tau$.
Rational functions $1/(z-v)$ with distinct $v$ are linearly independent,
as is seen by multiplying by $z-v$ and evaluating at $v$. Hence equality
of responses is exactly equality under $\mathcal L$, proving
\eqref{eq:rootcubefiber}. The left kernel of $B$ consists of vectors
constant on each connected component, so
$\rank\mathcal L=\rank B=|V(\mathcal P)|-\kappa(\mathcal P)$.
Its kernel vanishes precisely for a forest. If the kernel is nonzero,
a sufficiently small nonzero kernel vector and the zero vector lie in
the cube; \eqref{eq:petalrealization} realizes an actual collision.
At an interior point a neighborhood within the affine kernel remains in
the cube, giving the dimension formula. At boundary points the cube can
reduce the dimension, and the exact statement is the intersection
\eqref{eq:rootcubefiber}.
\end{proof}

In a forest, a leaf residue determines its unique incident $\tau_j$;
subtracting its contribution and deleting the edge gives recursive
recovery. Conversely, a cycle flow divided by the positive factors $s_j$
gives a kernel vector. Every loopless pole multigraph without isolated vertices, with distinct
positive labels on its vertices, is realizable: for an edge with endpoint values
$l<u$, choose $\alpha=\beta=l$ and $\gamma=(u-l)/2$.
Even when all $\tau_j$ are recovered, relative axes of different petal
holonomies can remain invisible. This is the continuous gauge ambiguity
already illustrated by cyclic blocks in Theorem~\ref{thm:fibers}.

\begin{example}[Three distinct denominators with an actual collision]
\label{ex:threepetals}
Take the three weight triples
\[
 (\alpha,\beta,\gamma)=(1,1,\tfrac12),\quad
 (2,2,\tfrac12),\quad(1,1,1).
\]
Then
\[
 H_1=(z-1)(z-2),\quad H_2=(z-2)(z-3),\quad H_3=(z-1)(z-3),
 \qquad(c_1,c_2,c_3)=(\tfrac12,2,1).
\]

\begin{figure}[htbp]
\centering
\begin{tikzpicture}[x=.95cm,y=.95cm,font=\small,line width=.65pt,>=stealth,
  vertex/.style={circle,fill=black,inner sep=1.7pt,outer sep=0pt}]
  \node[vertex,label=above:$0$] (o) at (0,0) {};
  \foreach \n/\x/\y in {a/-1.7/1.3,b/-1.8/-.2,c/1.7/1.3,
                          d/1.8/-.2,e/-.75/-1.7,f/.75/-1.7}
    \node[vertex] (\n) at (\x,\y) {};
  \draw (o)--node[above right] {$1$}(a)
    --node[left] {$\tfrac12$}(b)--node[below] {$1$}(o);
  \draw (o)--node[above left] {$2$}(c)
    --node[right] {$\tfrac12$}(d)--node[below] {$2$}(o);
  \draw (o)--node[left] {$1$}(e)
    --node[below] {$1$}(f)--node[right] {$1$}(o);
  \node at (-1.15,.5) {$j=1$};
  \node at (1.15,.5) {$j=2$};
  \node at (0,-1.2) {$j=3$};
  \node at (0,-2.35) {(a) $J$};
  \begin{scope}[xshift=6cm]
    \node[vertex,label=above:$1$] (p1) at (0,1.1) {};
    \node[vertex,label=below left:$2$] (p2) at (-1.15,-.85) {};
    \node[vertex,label=below right:$3$] (p3) at (1.15,-.85) {};
    \draw[->] (p1)--node[left] {$j=1$}(p2);
    \draw[->] (p2)--node[below] {$j=2$}(p3);
    \draw[->] (p1)--node[right] {$j=3$}(p3);
    \node at (0,-2.35) {(b) $\mathcal P$};
  \end{scope}
\end{tikzpicture}
\caption{The three-petal attachment and its pole graph in Example~\ref{ex:threepetals}. (a) The attachment $J$ rooted at $0$; edge labels are fixed weights, and $j$ indexes the petals. (b) Vertices of $\mathcal P$ are the numerical pole values of the denominators $H_j$; edge $j$ joins the two roots of $H_j$, directed from the lower pole to the higher one. The three simple cycles in $J$ correspond to the three edges of $\mathcal P$; a signed flow around an underlying undirected cycle of $\mathcal P$ produces a kernel direction of the response.}
\label{fig:poles}
\end{figure}
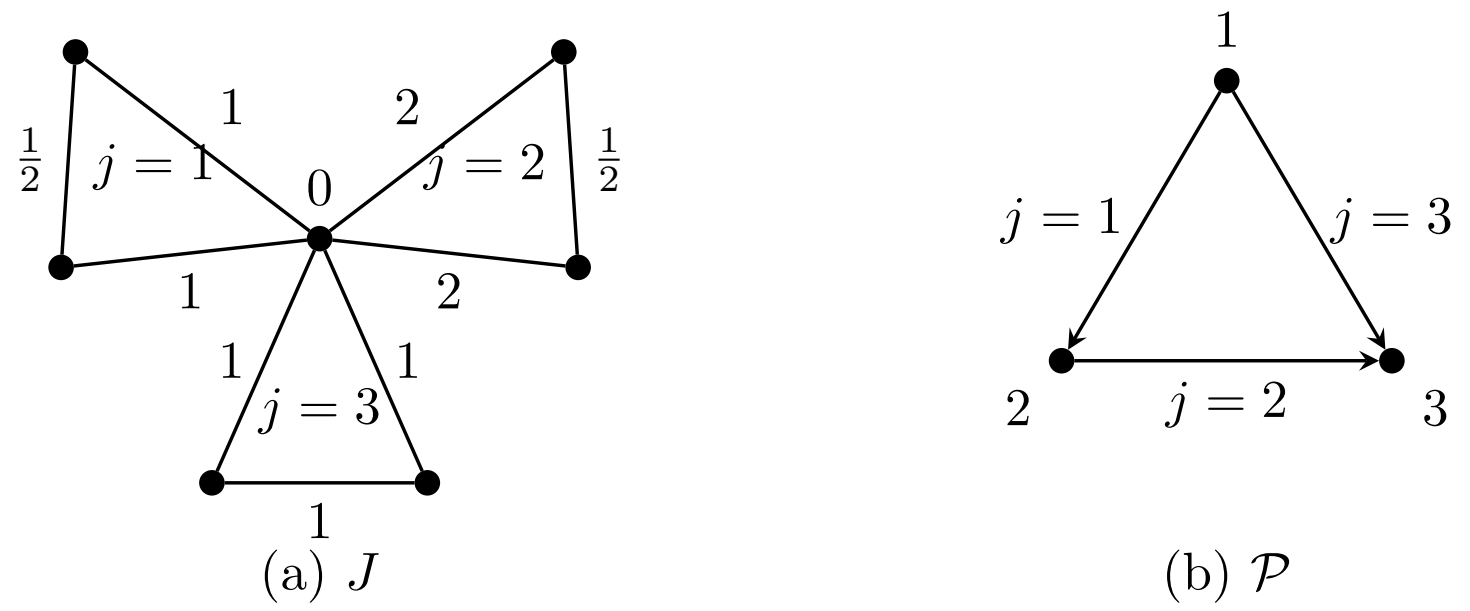

Figure~\ref{fig:poles} matches the attachment to the denominator poles.
The pole graph is a triangle, and
\[
 \frac1{H_1}+\frac1{H_2}-\frac2{H_3}
 =\frac{(z-3)+(z-1)-2(z-2)}{(z-1)(z-2)(z-3)}=0.
\]
Thus $\tau=(0,0,0)$ and $\tau'=(1/2,1/8,-1/2)$ have the same response;
\eqref{eq:petalrealization} supplies actual connections with these data.
Their different simple-cycle traces preclude gauge equivalence.
For both connections $R_\tau=0$, $\delta=8$, and
$(M_1,M_2,M_3)=(2z-3,8z-20,2z-4)$. Substitution in
\eqref{eq:petalresponse} gives the common reduced polynomial
\begin{equation}\label{eq:threepetalpoly}
 P_J=(z-1)(z-2)(z-3)
       (z^4-14z^3+47z^2-49z+11).
\end{equation}
The full characteristic polynomial of the $14$-dimensional complex
operator is $P_J^2$. By \eqref{eq:rootgluing}, the two attachments remain
isospectral after gluing to any identical core, for all allowed core
controls. Pairwise distinct denominators therefore do not suffice for
separation. The example fixes the attachment weights, and does not
contradict the all-edge-weight equivalence.
\end{example}

\subsubsection{Two-control core separation}
The converse passage from controlled spectra to attachment responses
uses the highest spectral coefficients in the two control variables.
These coefficients first identify four core traces, after which the
remaining comparison separates the attachment response.

Let the core be $K_4$ on $0,1,2,3$, rooted at $0$, with controls
$x=w_{01},u=w_{02}$ and fixed positive weights
$(w_{03},w_{12},w_{23},w_{31})=(v,p,q,h)$. Gauge its root links to the
identity and set $a=U_{12}$, $b=U_{23}$, $d=U_{31}$; in the transport
convention of Section~2, $a$ transports from $2$ to $1$.
The seven simple-cycle halftraces are
\[
 \operatorname{Re}a,\ \operatorname{Re}b,\ \operatorname{Re}d,
 \ \operatorname{Re}(ab),\ \operatorname{Re}(bd),
 \ \operatorname{Re}(da),\ \operatorname{Re}(abd).
\]
A fixed rooted attachment $J$ meets the core only at $0$, has $m$ nonroot
vertices, and has root degree $\delta>0$. Competitors use the same graph
and positive weights. Write $P_i,R_i$ for the bare core polynomials and
$D_i,C_i$ for the attachment pair, for $i=U,V$. For every such quantity
$Z$, write $\Delta Z=Z_U-Z_V$.

\begin{lemma}[Spectral extraction]\label{lem:coreextraction}
For a pair with equal spectra for every $x,u>0$, the four
quantities
\[
 \operatorname{Re}a,\quad\operatorname{Re}b,\quad
 \operatorname{Re}d,\quad\operatorname{Re}(bd)
\]
are common to the pair. With
\[
 d_3=v+h+q,\quad t=z-d_3,\quad k=d_3+2p(1-\operatorname{Re}a),
 \quad N=|v+hd^*+qb|^2,\quad Q=(3z-k)t-N,
\]
the common $k,N$ satisfy
\begin{equation}\label{eq:coreextraction}
 t\Delta C=-Q\Delta D,\qquad \deg\Delta D\le m-3.
\end{equation}
The degree bound includes the zero polynomial.
\end{lemma}
The proof is given in Appendix~\ref{app:core}.

\begin{proposition}[Response separation]\label{prop:coreconditional}
For an isospectral pair as above,
if
\begin{equation}\label{eq:coreweights}
 h\ne q\qquad\text{or}\qquad h=q\ \text{ and }\ v\ge h+2p,
\end{equation}
then $(D_U,C_U)=(D_V,C_V)$ and all seven core simple-cycle halftraces agree.
Without \eqref{eq:coreweights}, unequal responses can occur only if
\begin{equation}\label{eq:coresymmetric}
 h=q,\qquad \delta=h+p(1+\operatorname{Re}a)-v.
\end{equation}
This necessary condition does not assert realizability of unequal responses.
\end{proposition}
The proof is given in Appendix~\ref{app:core}.

For this core family, \eqref{eq:coreweights} reduces recovery of all
simple-cycle Wilson data to separation by the attachment map
$U_J\mapsto(D_J,C_J)$. Every simple cycle lies in one of the two pieces;
the core traces and the attachment response determine the controlled
polynomial by \eqref{eq:rootgluing}. Conversely, an actual response
collision glued to one fixed core gives an actual Wilson collision.
For triangular attachments, Proposition~\ref{prop:poleforest} therefore
gives a forest criterion for this controlled experiment.
Example~\ref{ex:threepetals} supplies a failure with actual connections.
Outside \eqref{eq:coreweights}, the necessary compensation condition
\eqref{eq:coresymmetric} can be tested by the general decision procedure
in Section~\ref{sec:arbitrarycontrols}.

\section{Scalar on-site terms and higher-rank limitations}\label{sec:extensions}
Let $v^0\in\R^n$ and $B\in\R^{n\times m}$ be known, and let $A_w(U)$
be the block adjacency matrix with $(x,y)$ block $w_{xy}U_{xy}$ on an
edge and zero otherwise. Consider
\begin{equation}\label{eq:onsite}
 H_w(U)=\operatorname{diag}(v^0+Bw)\otimes I_2-A_w(U).
\end{equation}
Here $\otimes$ is the Kronecker product, so the scalar potential at each
vertex acts on its two-dimensional fiber as a multiple of $I_2$.
The same diagonal data are used for all connections. The case $B=0$
permits fixed scalar on-site potentials, including the negative adjacency
operator when $v^0=0$.

Independent scalar diagonal variables encode principal minors through a
characteristic determinant. The image of the principal-minor map for
complex Hermitian matrices is studied in \cite{AlAhmadiehVinzant}.
In \eqref{eq:onsite}, the two coordinates at each vertex share the same
scalar potential.

\begin{proposition}\label{prop:onsite}
The all-weight spectral fibers of \eqref{eq:onsite} are exactly those
of $W_G$. The structural, full-fiber, finite-measurement and sparse-cycle
conclusions above remain valid. Stability holds with the corresponding
spectral norm bound. Effective finite-bit algebraic reconstruction applies
when $v^0,B$ are real algebraic.
\end{proposition}
\begin{proof}
The isolated-cycle moment $\Tr H_{w_C(t)}(U)^\ell$ still has degree at
most $\ell$ in $t$. At zero only the cycle remains off the diagonal.
Diagonal steps cannot occur in a walk winding once in $\ell$ steps, and
all nonwinding contributions are independent of the connection. Thus
its difference from the identity connection is exactly
\eqref{eq:isolation}. The reverse implication follows from the self-dual
determinant expansion with arbitrary scalar diagonal.

Moment degrees and quotient dimension are unchanged, so the coefficient
and secant proof transfers. On weights in $(0,1]$, all eigenvalues have
absolute value at most
\[
 R_* =\max_x\left(|v_x^0|+\sum_e|B_{xe}|\right)+\Delta.
\]
Replace $2\Delta+1$ by $R_*+1$ in \eqref{eq:noisechain}; the rest of
reconstruction uses the same cycle data. Exact algebraic decisions
require effective algebraic coefficients, as stated.
\end{proof}

For $SU(r)=\{Q\in\C^{r\times r}:Q^*Q=I_r,\det Q=1\}$, use the analogous Hermitian operator with scalar real
diagonal and scalar real edge weights, and let $\overline U$ be entrywise
complex conjugation of the transports.

\begin{proposition}\label{prop:higher}
For $r\ge3$, every graph containing a cycle has generically at least two
inequivalent connections with identical all-weight spectra. Forests are
exactly the graphs with generic uniqueness in this model. On a single
cycle, the generic entire fiber has exactly two classes.
\end{proposition}
\begin{proof}
Since $H_w(\overline U)=\overline{H_w(U)}$ and the matrix is Hermitian,
its real spectrum is unchanged. On a chosen cycle the condition
$\operatorname{Im}\operatorname{tr}\Hol_U(C)\ne0$ is generic: it is a
nonzero real-analytic condition, as witnessed by placing
$\operatorname{diag}(e^{ia},e^{ib},e^{-i(a+b)},1,\ldots,1)$ on one edge
for suitable $a,b$. The trace then differs from its conjugate, excluding
gauge equivalence. Forests have no free chord and always have a single class.

For a cycle of length $\ell$, tree gauge and diagonalization of the
holonomy give scalar magnetic cycles with phases
$\theta_1,\ldots,\theta_r$, whose sum is zero modulo $2\pi$.
For any fixed positive weights each scalar characteristic polynomial is
\[
 P_w(\lambda)+c_w\cos\theta_j,\qquad
 c_w=2(-1)^{\ell-1}\prod_{e\in C}w_e\ne0,
\]
where $P_w$ is known, monic and phase independent. Indeed, in
$\det(\lambda I-H)$ the two directed permutation cycles have sign
$(-1)^{\ell-1}$ and positive scalar weight factors multiplying the phase factors.
The full polynomial is $\prod_j(P_w(\lambda)+c_w\cos\theta_j)$ and
determines the multiset of cosines, since substitution of a nonconstant
polynomial is injective.

Generically these cosines are distinct and no nonempty proper subset
of holonomy eigenvalues has product $1$ or $-1$. One precise exceptional
set is given by the vanishing of the discriminant of $(Q+Q^*)/2$ or of
$\det(\bigwedge^k Q-I)\det(\bigwedge^k Q+I)$ for $1\le k<r$, where $Q$
is the holonomy and $\bigwedge^k Q$ is its induced action on the $k$th
exterior power of $\C^r$. These real-analytic conditions are proper: on the
determinant-one torus they exclude finitely many proper phase relations.
Their complement is open dense and full Haar measure.

Any competitor, after permutation, has eigenphases
$\varepsilon_j\theta_j$ with $\varepsilon_j\in\{1,-1\}$. Its
determinant-one condition gives
$1=\prod_j e^{i\varepsilon_j\theta_j}
=(\prod_{\varepsilon_j=-1}e^{i\theta_j})^{-2}$,
since $\prod_j e^{i\theta_j}=1$. Thus the product of the eigenvalues whose
signs were flipped lies in $\{1,-1\}$. The excluded proper-subset phase
relations therefore force either none or all to be flipped. On the additional generic nonreal-trace set these are two
distinct conjugacy classes, both realized and isospectral for all weights.
\end{proof}

The obstruction is the elementary antiunitary symmetry of a Hermitian
spectral problem. For $SU(2)$ complex conjugation is inner and produces
no new gauge class. The higher-rank fiber can be larger at special points:
for $SU(4)$ phases $(\pi/6,5\pi/6,\pi/4,3\pi/4)$, each indicated pair
has product $-1$. Reversing either pair independently gives four distinct
classes with identical all-weight cycle spectra.

\section{Further questions}
The preceding constructions give all exceptional fibers and an algorithmic
criterion for arbitrary restricted controls. Quantitative questions remain
about efficient implementation and graph descriptions of the resulting
semialgebraic conditions. In higher rank, rigidity should first
be considered modulo complex conjugation; additional ambiguities may
remain, and the quaternionic simple-cycle determinant reduction no
longer applies in the same form.

For restricted scalar controls, a nonzero
moment polynomial may vanish identically on a lower-dimensional set of
accessible weights. The finite-measurement argument then requires an
additional separation hypothesis. A direct graph criterion for that hypothesis
would address experiments with only a few tunable couplings. Quantitative
questions include the gap between $3\beta-3$ necessary scalar observations
and $4\beta-5$ selected cycles, and the conditioning of their reconstruction.

\section*{Acknowledgments}
The author thanks Professor Bobo Hua for his support.

\section*{Use of AI tools}
The author used OpenAI's Codex to investigate candidate constructions,
check calculations, and solicit separate AI reviews of draft proofs. The
author determined the direction of the work, verified all mathematical
claims and arguments, prepared the final manuscript, and assumes full
responsibility for its contents.

\bibliographystyle{amsplain}
\bibliography{references}
\appendix
\section{Certified reconstruction estimates}\label{app:stable}
\begin{proof}[Proof of Theorem~\ref{thm:noisy}]
For a tree output its unique gauge class with $E_G=0$. Otherwise we
construct local representatives in compatible tree gauges and then combine blocks.

\emph{$K_4$ and doubled-triangle anchors.}
For $K_4$ require $p=\det(\operatorname{Im}A,\operatorname{Im}B,
\operatorname{Im}C)\ne0$ and use the positive frame from its first two
imaginary vectors. For $D$ require $p=F_1F_2F_3\ne0$ and use the frame
in Lemma~\ref{lem:double}. Both canonical representatives are unique,
analytic and semialgebraic on these domains. Their seven- or ten-coordinate
observation map $h$ is injective and immersive. For $K_4$, positive
Gram reconstruction gives the differential rank. For $D$, nonangular
coordinates follow smoothly from the first seven data. Replacing the last three observations by the corresponding cosine values
(subtracting the known nonangular terms and dividing by the nonzero
transverse-length factors), the angular Jacobian has rows
$(-s_1,0),(0,-s_2),(s_3,-s_3)$, where
\[
 s_1=\sin(\theta_B-\phi),\quad s_2=\sin(\theta_C-\phi),\quad
 s_3=\sin(\theta_C-\theta_B-\phi).
\]
Rank failure forces two of these to vanish. The three cases force,
respectively, $\sin(\theta_C-\theta_B)=0$,
$\sin(\theta_C-2\phi)=0$, or $\sin\theta_B=0$, each excluded by $p\ne0$.

Let $y$ be the vector of measured anchor traces in the same coordinate
order as the observation map $h$. Consider the compact feasible set of unit
anchor tuples satisfying
$\|h(x)-y\|_\infty\le\eta$. It contains the true tuple. Effective real
quantifier elimination and sampling \cite[Section 2]{Basu} decide whether
this set avoids $p=0$. If not, report uncertified. Otherwise search
$\delta=1,1/2,1/4,\ldots$ for the first threshold with $|p|\ge\delta$
throughout the feasible set, and sample a feasible tuple. Positive square
roots in canonicalization are semialgebraic constraints.
Let $A_\delta$ be the set of canonical unit anchor tuples $x$ satisfying
$|p(x)|\ge\delta$. It is compact, since $p$ bounds its frame denominators
away from zero. Local inverses of the injective immersion
$h$ and compact separation of far pairs give a finite $C_\delta$ such that
\begin{equation}\label{eq:anchorbound}
 \max_j|x_j-x'_j|\le C_\delta\|h(x)-h(x')\|_\infty,
 \qquad x,x'\in A_\delta.
\end{equation}
Search $C=1,2,4,\ldots$ and test this universal semialgebraic inequality.
Existence ensures finite termination. Thus a sampled anchor has variable
error at most $2C_\delta\eta$ in the canonical gauge. For $K_4$ these
variables are the tree-gauge chords. For $D$ the actual chords are
$D_0,A,B,D_0C$, so the last chord error is at most the sum of the errors
in $D_0$ and $C$. The initial chord bound for the ear recursion is therefore
\[
 E_0=\begin{cases}2C_\delta\eta,&K_4\text{ anchor},\\
                  4C_\delta\eta,&D\text{ anchor}.
       \end{cases}
\]

\emph{Ears with uncertain old paths.}
Include the branch forest in the anchor tree and all but the last edge
of each new ear in the growing tree. Only its total transport $q_i$ is
new. In this proof $\|\cdot\|_2$ denotes the matrix operator norm induced by
the Euclidean vector norm. Let $R_i$ consist of the four true old path rows, and
$\widehat R_i$ the rows computed from recovered unit edges. If the old
chord error is at most $E_{i-1}$ and $L_i$ is the largest number of old
nontree edges on a selected path, telescoping products gives
\[
 \|R_i-\widehat R_i\|_2\le2L_iE_{i-1}.
\]
First test whether $\widehat R_i$ is invertible; otherwise report
uncertified. For its smallest singular value $\sigma_i$, choose a
certified bound $\sigma_i/2\le s_i\le\sigma_i$. Let $y_i$ be the vector
of the four measured cycle traces for this ear, in the row order of
$\widehat R_i$. Solve
$z_i=\widehat R_i^{-1}y_i$ and put
$h_i=(2\eta+2L_iE_{i-1})/s_i$. If $h_i\ge1$, report uncertified.
Otherwise, writing $c_i=R_iq_i$ for the true four traces, we have
\[
 \widehat R_i(z_i-q_i)=(y_i-c_i)+(R_i-\widehat R_i)q_i,
 \qquad |y_i-c_i|\le2\eta,
\]
so $|z_i-q_i|\le h_i<1$. In particular $z_i\ne0$, and
\[
 \left|\frac{z_i}{|z_i|}-q_i\right|
 \le |1-|z_i||+|z_i-q_i|\le2h_i.
\]
Use the recursion
\begin{equation}\label{eq:earbound}
 E_i=\max\left\{E_{i-1},\frac{4(\eta+L_iE_{i-1})}{s_i}\right\}.
\end{equation}
This accounts for errors in old paths, not just the new measurements.
Moreover success implies $2L_iE_{i-1}<s_i$, so the true frame is also
invertible. Therefore certificate success itself justifies the local
singleton fiber used in the construction.

\emph{Gram blocks and orientation pairs.}
For a block $B$ containing neither $K_4$ nor $D$ as a minor, write
$\beta_B=|E(B)|-|V(B)|+1$ and let $N=\beta_B+1$, $k=\min(4,N)$. Use the first
$k$ paths as references. Form their measured Gram matrix $A_y$, with
diagonal one. Denoting the smallest eigenvalue by $\lambda_{\min}$, test
$\rho=\lambda_{\min}(A_y)-3\eta>0$.
Sample any feasible unit tuple for the selected observations; the true
tuple guarantees nonemptiness. Both its reference matrix and the true
one have smallest eigenvalue at least $\rho$, since each row of their
difference from $A_y$ has at most three errors of size $\eta$.
Write these matrices as $A,A'$ and the cross-reference columns as
$v_i,v'_i$. Their entry differences are at most $\delta_g=2\eta$.
For $S=A^{1/2}$ and $S'=(A')^{1/2}$, the Sylvester identity and integral
solution give
\[
 \|S-S'\|_2\le\frac{3\delta_g}{2\sqrt\rho},\qquad
 \|S^{-1}-(S')^{-1}\|_2\le\frac{3\delta_g}{2\rho^{3/2}}.
\]
Indeed $S(S-S')+(S-S')S'=A-A'$, and the integral solution has factors
$e^{-tS}$ and $e^{-tS'}$, each bounded by $e^{-t\sqrt\rho}$. The inverse
bound follows by multiplication by the two inverses; commutativity is
not required.

For $k=4$, use the rows of $S$ as reference transports and
$v_i^TS^{-1}$ for all other paths. For $k<4$ there are only $k$ paths;
embed the square root in $\R^4$ by adding zero coordinates.
The difference of two nonreference rows satisfies
\[
 \begin{split}
 |v_i^TS^{-1}-(v'_i)^T(S')^{-1}|
 &\le |v_i-v'_i|\,\|S^{-1}\|_2
       +|v'_i|\,\|S^{-1}-(S')^{-1}\|_2\\
 &\le\frac{2\delta_g}{\sqrt\rho}
       +\frac{3\delta_g}{\rho^{3/2}},
 \end{split}
\]
since $|v'_i|\le2$ and $|v_i-v'_i|\le2\delta_g$.
Thus corresponding rows differ by at most
\begin{equation}\label{eq:grambound}
 E_Q=\left(\frac2{\sqrt\rho}+\frac3{\rho^{3/2}}\right)2\eta.
\end{equation}
This also covers reference rows, as $\rho\le1$.
For $k=4$, include the tuple obtained by the same fixed orientation-reversing
orthogonal map on every row. It pairs the two true and fitted $SO(4)$
classes with this error. For $k<4$ there is only one class.

Realize the tuples recursively on the network: set ordinary prefixes and
suffixes to identity, implement parallel pieces independently, and put
each leaf path transport on its last edge. Nonidentity edges are some
$Q_i$. After a fixed tree gauge, each chord is a known word in these
variables and inverses. If $J_B$ is the maximal word length, telescoping
gives block error $E_B=J_BE_Q$.

\emph{Global fibers and small-noise success.}
For each combination of local candidates, allow independent block
rotations and divide by the common rotation. This is a compact set
$\widehat\F$ of valid gauge classes. For a block reconstructed from a rigid anchor,
define $E_B$ to be its final ear bound (or $E_0$ if there are no added ears).
The reported global bound is $E_G=\max_B E_B$.
Pair local true and fitted
representatives and use the same rotations on both sides. Every point
on either side has a point on the other within $\max_BE_B$, proving
both Hausdorff inequalities. 
The output need not be the exact fiber of a common noisy observation
vector; normalized ear fits need not retain residual at most $\eta$ in
every coordinate.

For a fixed generic true connection, whole-fiber anchor uniqueness and
compactness force all small-noise feasible anchors to satisfy a common
$|p|\ge\delta_0>0$. The first successful dyadic choice then ranges over
finitely many values, so its inverse constants stay bounded. All true
ear frames are invertible; induction shows that recovered frames converge
to them, their $s_i$ stay bounded below, and the recursion is $O(\eta)$.
Reference Gram eigenvalues also stay positive. There are finitely many
steps, so the algorithm succeeds for all small enough noise. The
nondegeneracy conditions are generic by the path and anchor lemmas and
the arbitrary Gram coordinates.
\end{proof}

\section{Coefficients for the two-control core}\label{app:core}
\begin{proof}[Proof of Lemma~\ref{lem:coreextraction}]
Here $P,R$ denote the bare core polynomials of one connection, with
subscripts $U,V$ when both are compared. For a polynomial $f$, the notation
$[x^a u^b]f$ means the coefficient of $x^a u^b$, retaining all other
variables; $[xu]f$ abbreviates $[x^1u^1]f$.
Put $A=\operatorname{Re}a$, $B=\operatorname{Re}b$,
$D_0=\operatorname{Re}d$, $T=\operatorname{Re}(abd)$,
$E=\operatorname{Re}(ab)$, $F=\operatorname{Re}(bd)$ and
$G=\operatorname{Re}(da)$.
Write
\[
 d_0=z-x-u-v,\quad d_1=z-x-p-h,\quad
 d_2=z-u-p-q,\quad t=z-v-h-q.
\]
The reduced permutation expansion gives
\begin{align*}
 P={}&d_0d_1d_2t-x^2d_2t-u^2d_1t-v^2d_1d_2
       -p^2d_0t-h^2d_0d_2-q^2d_0d_1\\
    &+x^2q^2+u^2h^2+v^2p^2
       +2(xupAt+uvqBd_1+xvhD_0d_2+pqhTd_0)\\
    &-2(xpqvE+xuqhF+uphvG),\\
 R={}&d_1d_2t-p^2t-h^2d_2-q^2d_1+2pqhT.
\end{align*}
In particular, $[xu]P=(3z-k)t-N$ and $[xu]R=t$, where
\[
 N=v^2+h^2+q^2+2vhD_0+2vqB+2hqF.
\]
For example, the mixed coefficient before collecting terms is
\[
 [xu]P=(3z-v-2p-h-q+2pA)t
       -v^2-h^2-q^2-2vhD_0-2vqB-2hqF,
\]
which is $Q=(3z-k)t-N$. The quadratic terms in each single
control cancel, so both reduced polynomials are affine in each control.

The leading three coefficients of $D_i$ are connection independent:
the first two reduced spectral moments of the principal matrix depend
only on its diagonal degrees and squared edge norms. Thus
$\deg\Delta D\le m-3$. The root Schur complement gives
\[
 C_i/D_i=-\delta-s_i(z),\qquad
 s_i(z)=z^{-1}\sum_{0j\in E(J)}w_{0j}^2+O(z^{-2}).
\]
Since $D_i$ is monic and its first three coefficients are common,
\[
 C_i=-\delta D_i-
 \left(\sum_{0j\in E(J)}w_{0j}^2\right)z^{m-1}+O(z^{m-2}).
\]
Thus $\deg\Delta C\le m-2$. Here and below a negative degree bound
means that the polynomial is zero, so attachments with $m<3$ are included. These are polynomial degree bounds obtained at
infinity and impose no condition at finite poles.

Equality of the complex characteristic polynomials implies equality of
their monic reduced square roots. The $xu$ coefficient of the gluing
identity therefore reads
\[
 0=\Delta D\,Q_U+D_V(2p\Delta A\,t-\Delta N)+t\Delta C.
\]
Its $z^{m+1}$ and then $z^m$ coefficients give $\Delta A=0$ and
$\Delta N=0$, proving \eqref{eq:coreextraction}.
In the full identity
\[
 0=\Delta D\,P_U+\Delta C\,R_U+D_V\Delta P+C_V\Delta R,
\]
the first two terms have $x$ and $u$ coefficients of degree at most $m$.
Also $\Delta R=2pqh\Delta T$ has no controls. The $xz^{m+1}$ and
$uz^{m+1}$ coefficients are therefore $2vh\Delta D_0$ and $2vq\Delta B$,
respectively.
Indeed, after $\Delta A=\Delta N=0$, the core expansion gives
\[
 [x^1u^0]\Delta P=2vh\Delta D_0\,z+O(1),\qquad
 [x^0u^1]\Delta P=2vq\Delta B\,z+O(1).
\]
Since $D_V$ is monic of degree $m$, these are the only contributions
in degree $m+1$; hence both vanish. The displayed formula for $N$,
and $hq>0$, then give
$\Delta F=0$. All four traces are common, including when $N=0$.
\end{proof}

\begin{proof}[Proof of Proposition~\ref{prop:coreconditional}]
Put $T=\Delta\operatorname{Re}(abd)$,
$E=\Delta\operatorname{Re}(ab)$, and $G=\Delta\operatorname{Re}(da)$.
The simple-cycle expansion of the core, using the four common traces, gives
\begin{align*}
 \Delta P&=2pqh(z-v)T-2xpq(hT+vE)-2uph(qT+vG),\\
 \Delta R&=2pqhT.
\end{align*}
Indeed, the remaining triangle $123$ contributes
$2pqh(z-x-u-v)T$ to $\Delta P$; the two remaining four-cycles contribute
$-2xpqvE$ and $-2uphvG$. Deleting the root leaves only the triangle term
$2pqhT$. Define $H=tP_U-QR_U$. Subtracting the two gluing formulas,
multiplying by $t$, and using \eqref{eq:coreextraction} gives
\begin{align}
 0={}&\Delta D H+2pqh\,t[D_V(z-v)+C_V]T\notag\\
    &-2xpq\,tD_V(hT+vE)-2uph\,tD_V(qT+vG).
 \label{eq:corefullidentity}
\end{align}
For $H_x=[x^1u^0]H$, $H_u=[x^0u^1]H$, and $H_0=[x^0u^0]H$,
the two single-control coefficients yield
\begin{equation}\label{eq:corelinear}
 \Delta D H_x=2pq(hT+vE)tD_V,\qquad
 \Delta D H_u=2ph(qT+vG)tD_V.
\end{equation}

We need only the first two coefficients of these polynomials. With
$s=v+p+q+h$ and $r_0=2p+h+q+d_3$, they follow from
\begin{align*}
 P_U|_{x=u=0}&=z^4-2sz^3+O(z^2),\\
 R_U|_{x=u=0}&=z^3-r_0z^2+O(z),\\
 [x^1u^0]P_U&=-2z^3+(4s-v-p-h)z^2+O(z),\\
 [x^0u^1]P_U&=-2z^3+(4s-v-p-q)z^2+O(z),\\
 [x^1u^0]R_U&=-z^2+(p+q+d_3)z+O(1),\\
 [x^0u^1]R_U&=-z^2+(p+h+d_3)z+O(1).
\end{align*}
Here normalized trace means one half of the complex matrix trace.
For example, the normalized trace of the bare core at $x=u=0$ is $2s$,
and the derivative in $x$ of the normalized trace of its square is
$2(v+p+h)$. The second elementary coefficient therefore has $x$
coefficient $4s-v-p-h$. In the root-deleted matrix, selecting the $-x$
diagonal entry leaves $-[(z-p-q)t-q^2]$, proving the corresponding formula
for $R_U$. Exchanging $1$ and $2$ gives the $u$ formulas. Substitution gives
\begin{align}
 H_x&=z^4-(v+h+3q+k)z^3+O(z^2),\notag\\
 H_u&=z^4-(v+3h+q+k)z^3+O(z^2),\label{eq:coreleading}\\
 H_0&=-2z^5+B_0z^4+O(z^3),\qquad
 B_0=4p+3h+3q+3d_3+k.\notag
\end{align}

Suppose $\Delta D$ is nonzero, with degree $\ell$ and leading coefficient
$c$. Since $H_x$ is monic, the first constant factor on the right of
\eqref{eq:corelinear} is nonzero. Degree comparison with monic $D_V$
gives $\ell=m-3$, while leading coefficients in both equations give
$2pq(hT+vE)=2ph(qT+vG)=c$. Thus $H_x=H_u$, and their cubic coefficients
force $h=q$. The constant-control part of \eqref{eq:corefullidentity}
has leading coefficient $-2c+2pqhT$ in degree $m+2$, so $T=c/(pqh)$.
Substituting $T=c/(pqh)$ and $\Delta D H_x=ctD_V$ into the
constant-control equation gives, in the rational function field,
\[
 0=\frac{ctD_VH_0}{H_x}
      +2ct[D_V(z-v)+C_V].
\]
Here $c\ne0$, and $t,D_V,H_x$ are nonzero polynomials, so cancellation
is legitimate and yields
\begin{equation}\label{eq:coreforcedresponse}
 \frac{C_V}{D_V}=-\frac{H_0}{2H_x}-(z-v).
\end{equation}
For $h=q$, set $\alpha=-(v+4h+k)$. Equation~\eqref{eq:coreleading} gives
$H_0/H_x=-2z+(B_0+2\alpha)+O(z^{-1})$. Hence the constant term on the
right of \eqref{eq:coreforcedresponse} is
\[
 v-B_0/2-\alpha=v-h-p(1+\operatorname{Re}a).
\]
Comparison with \eqref{eq:rootschur} proves \eqref{eq:coresymmetric}.
Since $\delta>0$ and $\operatorname{Re}a\le1$, condition
\eqref{eq:coreweights} excludes this possibility, including its equality
case. No division by $N$ was used.

It follows that $\Delta D=0$, and \eqref{eq:coreextraction} then gives
$\Delta C=0$. The unmultiplied gluing identity has constant-control
coefficient $2pqh[D_V(z-v)+C_V]T=0$. The bracket is monic of degree
$m+1$, so $T=0$; its $x$ and $u$ coefficients then give $E=G=0$.
Together with the four common traces this proves the assertion.
\end{proof}

\end{document}